\documentclass[10pt,a4paper]{article}

\usepackage{arxiv}
\usepackage{tikz}
\usepackage[utf8]{inputenc} 
\usepackage[T1]{fontenc}    
\usepackage{hyperref}       
\usepackage{url}            
\usepackage{booktabs}       
\usepackage{amsfonts}       
\usepackage{nicefrac}       
\usepackage{microtype}      
\usepackage{lipsum}

\usepackage{amssymb}
\usepackage{lineno,hyperref}
\usepackage{amsmath}
\usepackage{bm}
\usepackage{subcaption}
\usepackage{graphicx}
\usepackage{svg}
\usepackage{rotating}
\usepackage{lscape}
\usepackage{array,multirow}
\usepackage[title]{appendix}
\usepackage{float}
\usepackage{stmaryrd}
\usepackage{textcomp}
\usepackage{algorithm}
\usepackage{algpseudocode}
\usepackage{amsthm}
\usepackage{listings}
\usepackage{xcolor}
\usepackage{siunitx}

\usepackage{pifont}
\usepackage{lineno}

\usepackage{import}
\usepackage{physics}

\usetikzlibrary{positioning}

\usepackage[]{natbib}
\author{
  Jan N. Fuhg \\
  Sibley School of Mechanical and Aerospace Engineering \\
  Cornell University, 
   New York, USA \\
  \texttt{jf853@cornell.edu} \\
   \And
   Ioannis Kalogeris \\
  School of Civil Engineering \\
  National Technical University of Athens, Greece
   \And
   Am\'{e}lie Fau\\
   Universit\'{e} Paris-Saclay, ENS Paris-Saclay,\\ CNRS,  LMT, 
Laboratoire de Mécanique et Technologie, France
   \And
 Nikolaos Bouklas \\
  Sibley School of Mechanical and Aerospace Engineering\\
  Center for Applied Mathematics\\
  Cornell University,
   New York, USA 
  }

\title{Interval and fuzzy physics-informed neural networks for uncertain fields}

\DeclareMathOperator*{\argmin}{arg\,min}

\begin{document}
\maketitle
\begin{abstract}
Temporally and spatially dependent uncertain parameters are regularly encountered in engineering applications.
Commonly these uncertainties are accounted for using random fields and processes, which require knowledge about the appearing probability distributions functions that is not readily available.
In these cases non-probabilistic approaches such as interval analysis and fuzzy set theory are helpful uncertainty measures. Partial differential equations involving fuzzy and interval fields are traditionally solved using the finite element method where the input fields are sampled using some basis function expansion methods. This approach however is problematic, as it is reliant on knowledge about the spatial correlation fields. 
In this work we utilize physics-informed neural networks (PINNs) to solve interval and fuzzy partial differential equations. The resulting network structures termed interval physics-informed neural networks (iPINNs) and fuzzy physics-informed neural networks (fPINNs) show promising results for obtaining bounded solutions of equations involving spatially and/or temporally uncertain parameter fields. In contrast to finite element approaches, no correlation length specification of the input fields as well as no Monte-Carlo simulations are necessary. In fact, information about the input interval fields is obtained directly as a byproduct of the presented solution scheme.
Furthermore, all major advantages of PINNs are retained, i.e. meshfree nature of the scheme, and ease of inverse problem set-up.
\end{abstract}

\keywords{Physics-informed machine learning \and Fuzzy set theory \and Interval set theory \and Non-probabilistic uncertainty }



\section{Introduction}
Uncertain parameters in engineering applications are commonly found to have temporal and spatial dependencies.
These uncertainties include external loading, e.g. in the case of wind loads, as well as material properties e.g. as found in heterogeneous media and biological tissues.
Commonly, these types of uncertainties have been modeled as stochastic or random fields \citep{ostoja1998random, ostoja2006material}. In stochastic fields the randomness is characterized by correlated random variables at each time increment and spatial location. When random fields appear in partial differential equations (PDEs) they generalize to stochastic PDEs (sPDEs), which are commonly solved using the so-called stochastic finite element method
\citep{der1988stochastic, stefanou2009stochastic}.
In this context, the conceptually simplest approach for uncertainty propagation are Monte-Carlo simulations, which are still commonly applied as the main benchmarking approach for probabilistic methods. However they generally require a large number of samples, which typically leads to high computational costs. In order to combat this, surrogate modeling techniques are used to vastly decrease the necessary computation time \citep{fuhg2019adaptive,fuhg2020state}.

A major bottleneck of probabilistic methods is the requirement of precise knowledge of the probability distribution functions of the uncertain parameters. This generally necessitates a great amount of statistical information in form of experimental data which is hard to come by.
Furthermore, following \cite{moller2000fuzzy} and \cite{bothe2013fuzzy}, stochastic uncertainty can also only be used if an event, which is a random result of an experiment, can be observed under constant boundary conditions. However, if the boundary conditions are (seemingly) changing or a non-significant amount of data is available there is a crucial lack of information which is called informal uncertainty.
Additionally, in certain problems, available information is only known with lexical imprecision.
Informal and lexical uncertainties are commonly present in engineering science and are areas where classical probabilistic methods may only be applied to a limited extend \citep{moller2003safety}.

For this reason, non-probabilistic approaches such as
interval analysis \citep{moore1966interval, jaulin2001interval, moore2009introduction} and
fuzzy set theory \citep{ZADEH1965338, klir1995fuzzy, zimmermann2011fuzzy} have been proposed. They can be used when informal or lexical uncertainties are present by describing these uncertainties through interval and fuzzy parameters. 
PDEs involving interval or fuzzy parameters generalize to interval and fuzzy PDEs (iPDEs and fPDEs).
Similarly to sPDEs, these constructs can be solved numerically with variants of the finite element method, e.g. 
interval finite element analysis \citep{chen2000interval, sofi2016novel, ni2020interval} and fuzzy finite element analysis \citep{rao1995fuzzy, muhanna1999formulation, yin2016fuzzy}.
However, interval finite element and fuzzy finite element methods,  as they pertain to the analysis of spatial and time-variable uncertainty fields, are still areas of active research since they are accompanied by high computational complexity \citep{schietzold2019}.
For example, efficient and convenient constructions of the interval fields are still posing a challenge \citep{ni2020interval} since they are typically constructed following some basis function expansion approach, based on the Karhunen–Loève (KL) decomposition \citep{moens2011numerical,verhaeghe2013interval,muscolino2013one, sofi2015static, ni2020interval,dannert2021}. Efficient intrusive approaches for uncertainty propagation have been developed based on these series expansions \citep{ni2020interval,SOFI2019} but they are applicable to linear problems with affine dependence on the interval fields. An additional limitation of KL-like expansions is that they are restricted to homogeneous fields and, more generally, they are intrinsically dependent on knowledge of the spatial correlation fields of the input, which might be assumed but are in general (without experimental data) not known \citep{dannert2018}. 
In this work we present a method based on artificial intelligence that completely bypasses the issues of all existing approaches, and obtains a description of the input field as a byproduct of the solution procedure. 

Machine learning technologies have been an emergent tool in computational engineering in recent years. They have  been used for data-driven constitutive modeling
\citep{huang2020machine, fuhg2021modeldatadriven, fuhg2021local,fuhg2021physics} as a means to enable hierarchical multiscale calculations or for the development of intrusive and non-intrusive Reduced Order Modeling (ROM) schemes for accelerated solutions of PDEs \citep{kadeethum2021non,hernandez2021deep, kadeethum2021framework}.
Recently, machine learning has also been employed as a solution scheme for PDEs. This approach called
physics-informed neural networks (PINNs) \citep{raissi2019physics,jagtap2020adaptive,wang2021and,fuhg2021mixed} has emerged as an alternative numerical technique for solving PDEs in a forward and inverse manner. Here, in an unsupervised manner, neural networks are used to approximate PDE solutions utilizing automatic differentiation to define global shape functions for the solutions. This formulation has already been extended to sPDEs with uncertain fields in \cite{yang2019adversarial}, \cite{zhang2019quantifying} and \cite{yang2021b}. One bottleneck with these approaches is that they still require information about the probability distribution functions which in turn requires a statistically significant amount of experimental data.

In this work we present a method for solving iPDEs and fPDEs based on machine learning that completely bypasses the issues of all existing approaches and obtains a description of the input field as a byproduct of the solution procedure. Specifically, we introduce interval physics-informed neural networks (iPINN), which are also extended here to fuzzy physics-informed neural networks (fPINN). The proposed iPINN framework consists of a dedicated network architecture and an appropriately modified loss function in order to obtain the extrema of the PDE’s solution field, while at the same time satisfying the physics of the problem. Importantly, the proposed methodology can be applied to any type of problem (linear or nonlinear) and any type of interval field (homogeneous or nonhomogenous) and fuzzy fields with convex membership functions, without requiring any information about the correlation structure of these fields. Furthermore, its ease of implementation and generalization capabilities, allows it to be straightforwardly applied to any type of fuzzy and interval PDE encountered in engineering applications. In this frame, iPINN can be understood as
a component of fPINN, on which we will elaborate in an upcoming section.

The paper is organized as follows. Interval and fuzzy variables and fields are shortly reviewed in Section \ref{sec::1}. PDEs involving these types of uncertainty measures are then reviewed in Section \ref{sec::2}. The iPINN and fPINN formulations are introduced in Section \ref{sec::3}, subsequently, these formulations are applied and studied through a series of numerical examples in Section \ref{sec::4}.
The paper is concluded in Section \ref{sec::5}.

\section{Interval and fuzzy variables and interval and fuzzy fields}\label{sec::1}
Even though interval and fuzzy set analysis has seen interest since the 1960s \citep{ZADEH1965338}, they have not received the same attention across all fields of computational physics. For this reason we use this section to describe the core concepts of fuzzy and interval numbers. 
An uncertain parameter $\bm{X}$ can be represented as an interval parameter $\bm{X}^{I}$ if it can be defined by a variation range 
\begin{equation}
   \bm{X} \in \bm{X}^{I} = \left[ \bm{X}^{L}, \bm{X}^{U} \right] 
\end{equation}
with the lower and upper limits $\bm{X}^{L}$ and $\bm{X}^{U}$ respectively.
An extension to this idea are interval fields.
Consider a spatial location $\bm{x} \in \mathcal{D}$ in a domain $\mathcal{D} \subset \mathcal{R}^{d}$ and let the time be defined by $t\in \mathcal{T}$.
A time-dependent parameter field $\bm{P}(\bm{x},t)$ is a spatially uncertain interval field $\bm{P}^{I}(\bm{x},t)$ if for any $\bm{x} \in \mathcal{D}$ and $t\in \mathcal{T}$ the field values are exactly defined by a spatially and time-dependent interval
\begin{equation}
    \bm{P}(\bm{x},t) \in \bm{P}^{I}(\bm{x},t)  = \left[ \bm{P}^{L}(\bm{x},t), \bm{P}^{U}(\bm{x},t)   \right].
\end{equation}
\begin{figure}[ht]
    \begin{subfigure}[b]{0.5\linewidth}
        \includegraphics[scale=0.22]{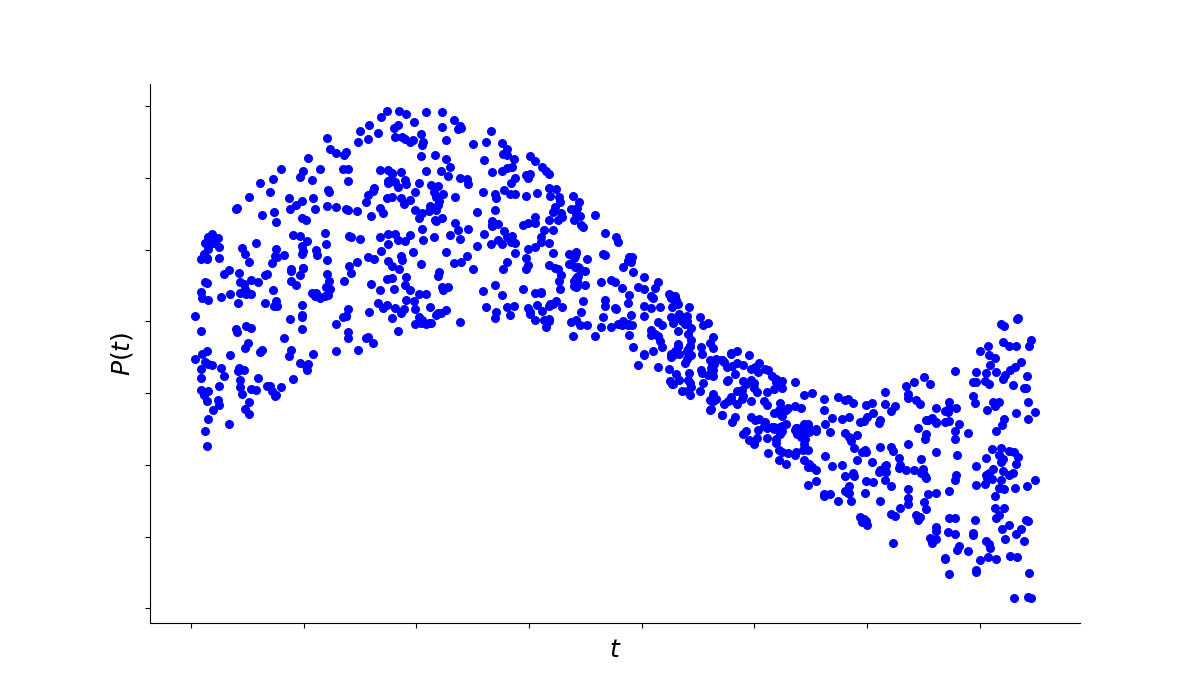}
    \caption{}\label{fig::IntFieldExp}
    \end{subfigure}
        \begin{subfigure}[b]{0.5\linewidth}
        \includegraphics[scale=0.22]{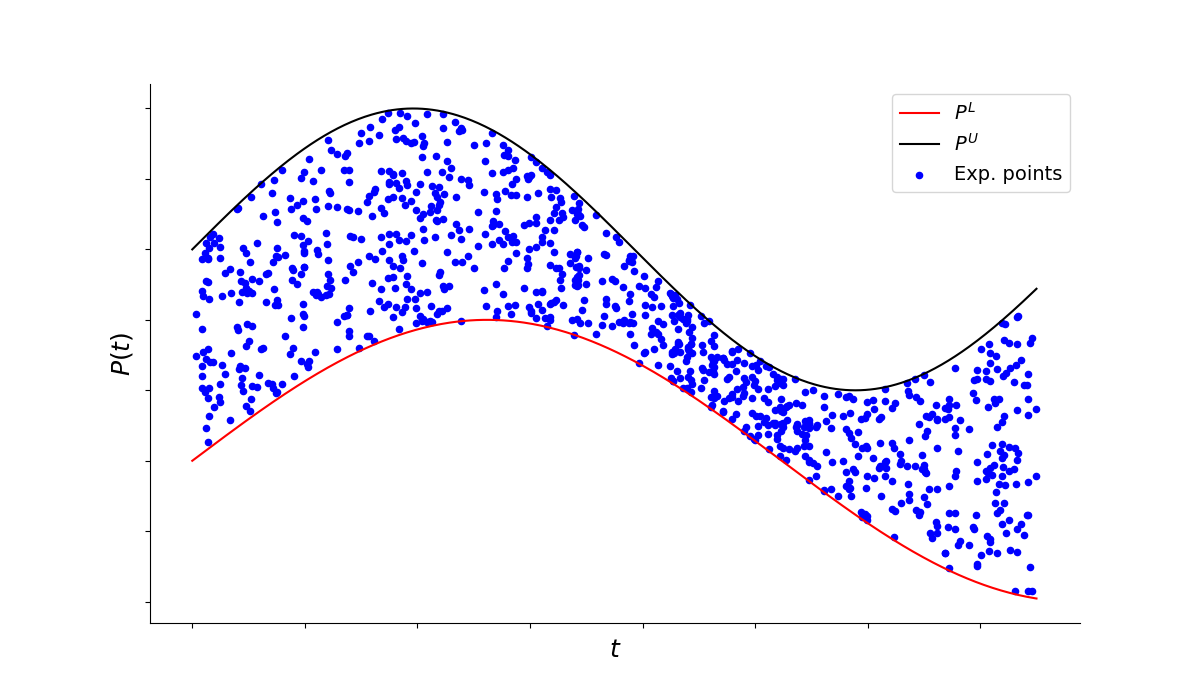}
    \caption{}\label{fig::IntFieldLines}
    \end{subfigure}
    \caption{Visual explanation of interval field concept. (a) experimental points, and (b) resulting approximated interval field.}
\label{}
\end{figure}
To make this more clear, assume we have measurements of a time-dependent process $P(t)$ (Figure \ref{fig::IntFieldExp}) where informal uncertainties are present (e.g. not enough data is available). To account for the uncertainty present in the data with interval theory we can construct bounding limits $P^{L}$ and $P^{U}$ (Figure \ref{fig::IntFieldLines}) which can for example be approximated using Gaussian process regression. Using these bounds we can assume restrictions on possible values of $P(t)$. This information can be used to obtain extreme values of variables that are dependent on $P(t)$. This is for example crucial for applications in engineering science, where structural failure can occur. 

From the provided definitions we can see that interval numbers can be understood as generalizations of real numbers. A further generalization are fuzzy numbers which can be seen as an extension to the interval concept. 
\begin{figure}[ht]
    \begin{subfigure}[b]{0.5\linewidth}
        \includegraphics[scale=0.75]{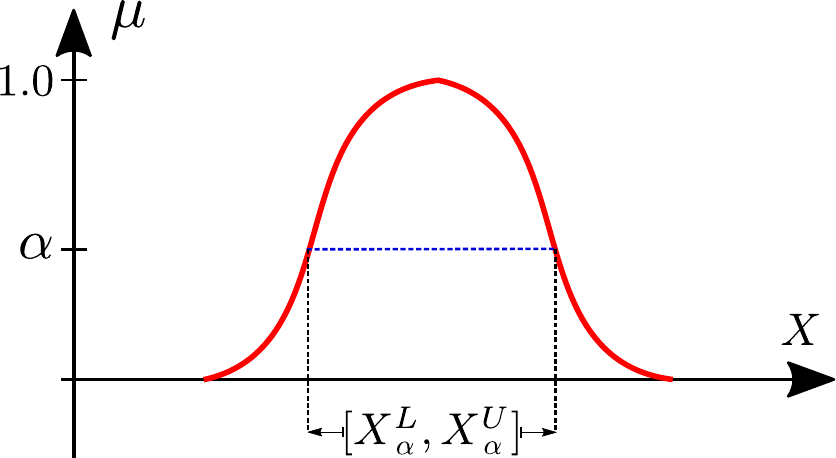}
    \caption{}\label{fig::FuzzyNumber}
    \end{subfigure}
        \begin{subfigure}[b]{0.5\linewidth}
        \includegraphics[scale=0.22]{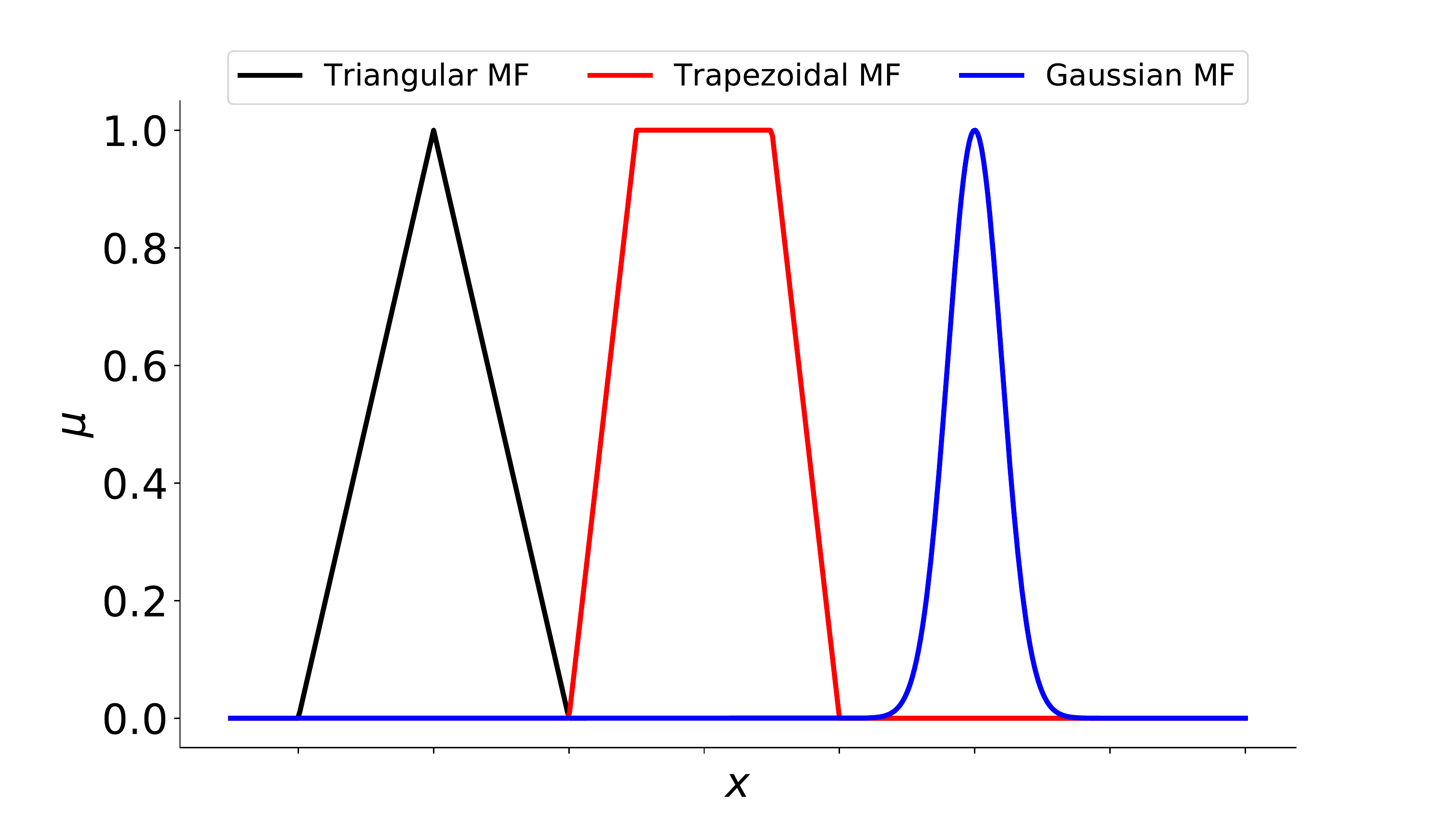}
    \caption{}\label{fig::FuzzyMemb}
    \end{subfigure}
    \caption{(a) Fuzzy number and $\alpha$-cut resulting in an interval variable, and (b) different types of fuzzy membership functions including triangular, trapezoidal and Gaussian.}
\label{}
\end{figure}

\begin{figure}[ht]
    \begin{subfigure}[b]{0.5\linewidth}
        \includegraphics[scale=0.22]{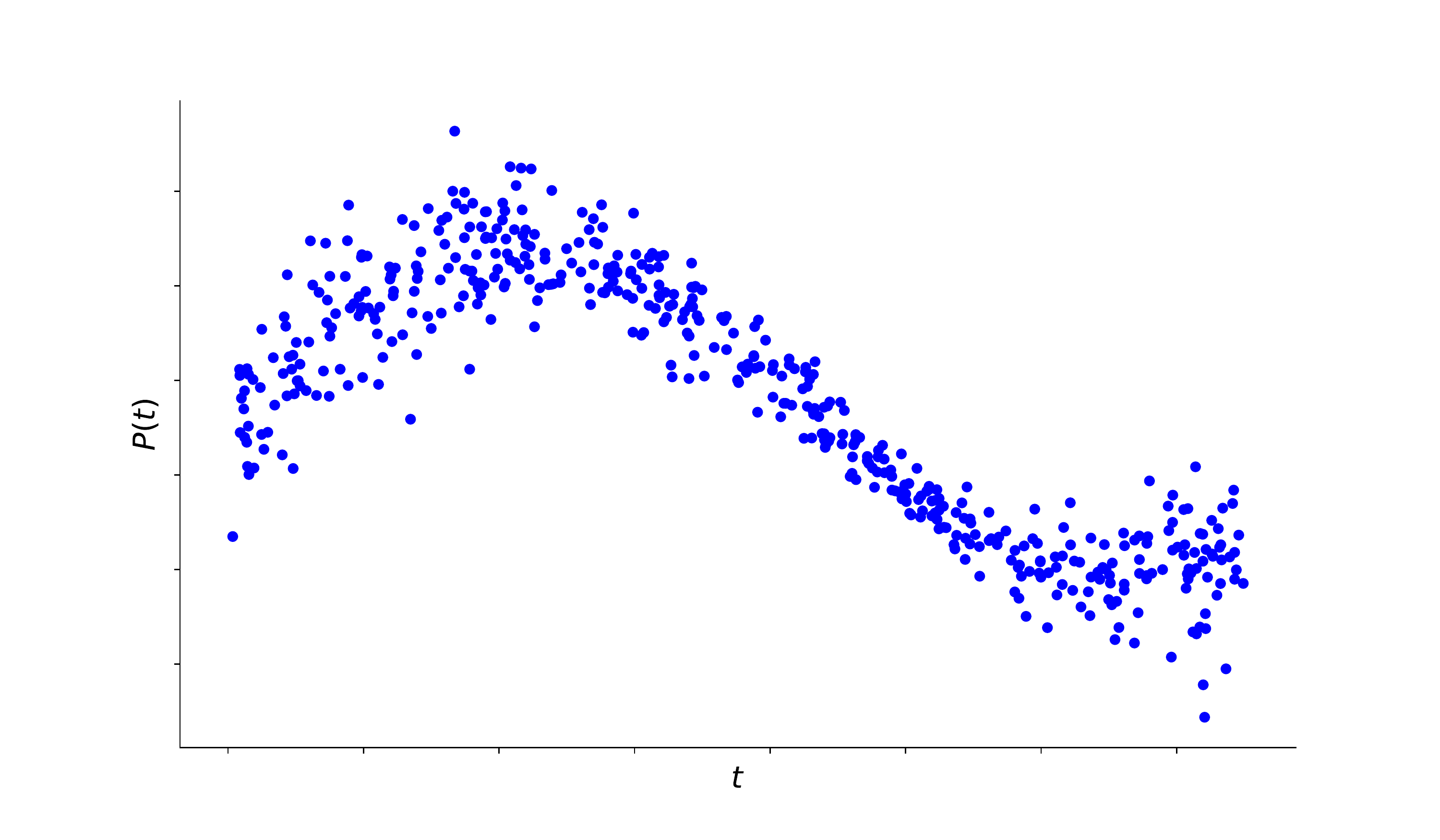}
    \caption{}\label{fig::fuzzyFieldExp}
    \end{subfigure}
        \begin{subfigure}[b]{0.5\linewidth}
        \includegraphics[scale=0.22]{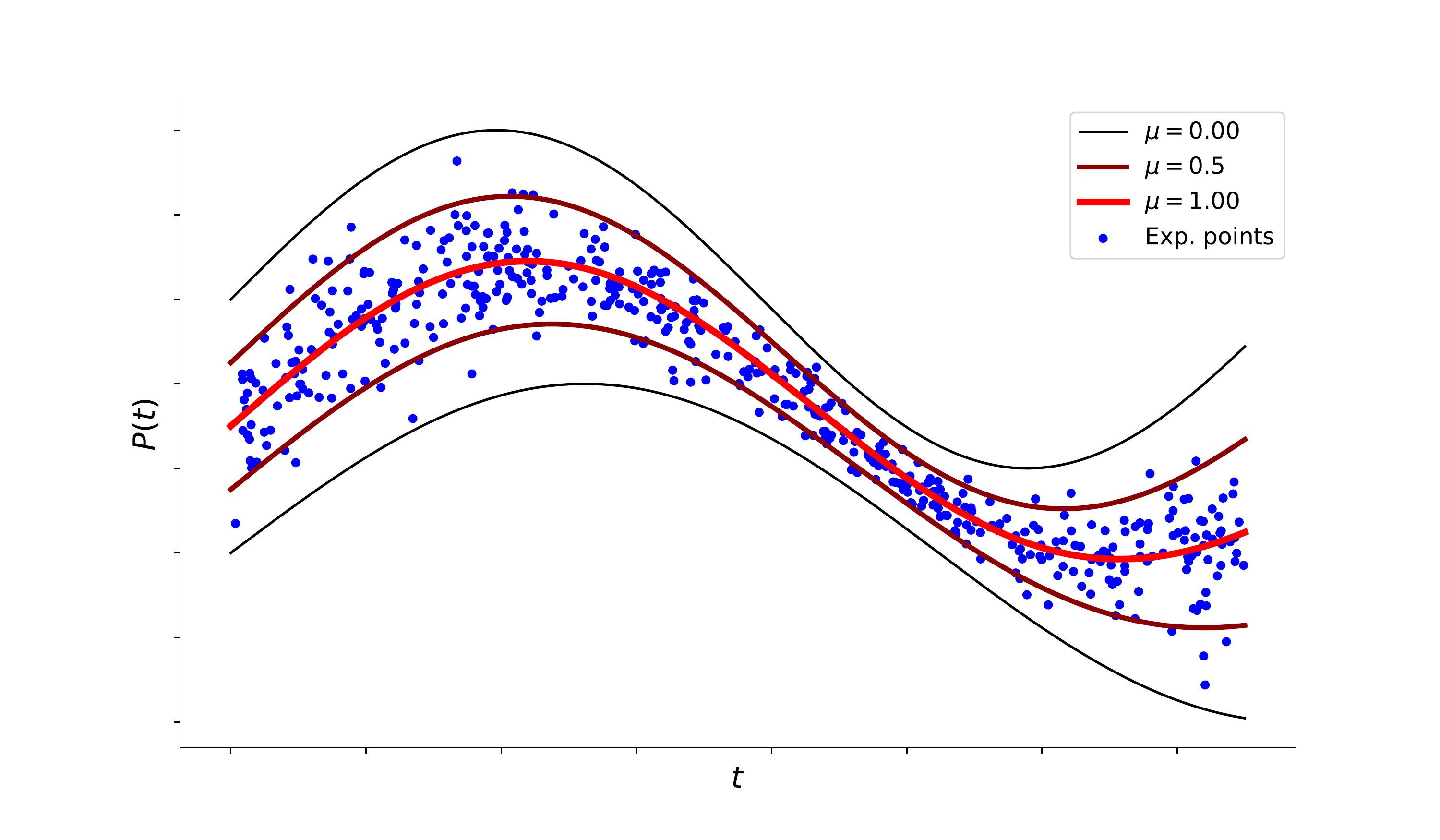}
    \caption{}\label{fig::fuzzyFieldLines}
    \end{subfigure}
    \caption{Visual explanation of fuzzy field concept. (a) experimental points, and (b) possible resulting membership contours of approximated fuzzy field.}
\label{fig::}
\end{figure}

A fuzzy number $\bm{X}^{F}$ is composed of the pair $\left( \bm{X}, \mu \right)$, with the membership function $\mu : \bm{X} \rightarrow [0,1]$. In contrast to the interval theory where an element is either part of a set or not, in fuzzy set theory the membership function $\mu$ allows us to assign a degree of membership to the elements of a set. This can be achieved by defining
\begin{equation}
    \begin{cases}
    \mu(\bm{X}) = 0, & \bm{X}\text{ is no member of the set}, \\
    \mu(\bm{X}) = 1, & \bm{X}\text{ is full member of the set}, \\
    0 < \mu(\bm{X}) < 1, & \bm{X}\text{ is partial member of the set}.
    \end{cases}
\end{equation}
This conceptually allows us to define a deterministic approach to account for uncertainties, for example by defining $\mu(\bm{X}_{1})= 0.2$ and $\mu(\bm{X}_{2})=0.8$ we can make an expert decision that $\bm{X}_{1}$ is less likely than $\bm{X}_{2}$.\\
The elements that are part of the set by at least a value $\alpha$ define the so-called $\alpha$-cut set 
\begin{equation}
    \bm{X}^{F}_{\alpha} = \lbrace \bm{x} \in \bm{X} | \mu(\bm{X}) \geq \alpha \rbrace.
\end{equation}
If the membership function is assumed to be convex, the $\alpha$-cut set is exactly defined by an interval variable
\begin{equation}
    \bm{X}^{I}_{\alpha} = \left[ \bm{X}_{\alpha}^{L}, \bm{X}_{\alpha}^{U} \right] ,
\end{equation}
see Figure \ref{fig::FuzzyNumber}. Different forms of convex membership functions can be found in the literature \citep{ZADEH1965338} including triangular, trapezoidal and Gaussian (Figure \ref{fig::FuzzyMemb}).
A fuzzy number can be thought of as exactly defined by an infinite number of $\alpha$-cuts, whereas a finite number of cuts $\bm{X}_{\alpha^{k}}$, $k=1,\ldots,r$ yields an approximation of $\bm{X}^{F}$. This concept allows us to obtain information of fuzzy uncertainties by simply relying on interval theory.\\
For example we can obtain an approximation of the fuzzy output of an operation $\bm{Y}^{F} = \bm{f}(\bm{A}^{F},\bm{B}^{F})$ between two fuzzy uncertain numbers by using $\alpha$-level optimization \citep{moller2000fuzzy}.
Assume that each element $i=1,\ldots,d$ in $\bm{A}^{F}$ and $\bm{B}^{F}$ has the same $\alpha$ level representation $A_{i,\alpha^{k}}$ and $B_{i,\alpha^{k}}$ where $k=1,\ldots,r$. The resulting $\alpha$-cuts of the output $Y_{i,\alpha^{k}}$ are then exactly defined by its smallest and largest element
\begin{equation}
    \begin{aligned}
        &\bm{Y}_{\alpha^{k}} = &&[\bm{Y}_{\alpha^{k}}^{L}, \bm{Y}_{\alpha^{k}}^{U}], \\
        &\text{with}&&
        Y_{i,\alpha^{k}}^{L} = \min_{\bm{a}\in \bm{A}_{\alpha^{k}}, \bm{b} \in \bm{B}_{\alpha^{k}}} f_{i}(\bm{a},\bm{b}), \\
      &  &&Y_{i,\alpha^{k}}^{U} = \max_{\bm{a}\in \bm{A}_{\alpha^{k}}, \bm{b} \in \bm{B}_{\alpha^{k}}} f_{i}(\bm{a},\bm{b}).
    \end{aligned}
\end{equation}
In an equivalent way to interval fields we are also able to define time and spatially dependent fuzzy fields as the pair $(\bm{F}(\bm{x},t), \mu)$ where $\mu : \bm{F}(\bm{x},t) \rightarrow [0,1]$. 
In order to understand how this allows us to account for informal and lexical uncertainties assume that we have access to experimental data of another time dependent field $P(t)$ for which we do not have enough statistical data to correctly approximate its moments. However, we have access to a rough approximation of the mean and user experiences such as lexical information (e.g. "Most of the time $P(t)$ behaves like ...") which allows us to account for the uncertainties in the dataset through a fuzzy field approach, see Figure \ref{fig::fuzzyFieldLines}. For more information on fuzzy set theory we refer to \cite{dubois2000fuzzy}.

\clearpage
\section{Interval and fuzzy partial differential equations}\label{sec::2}
Consider a fuzzy PDE of the form
\begin{equation}
    \tilde{\bm{G}}(\bm{x},t, \bm{u}, \mathcal{F}) = \bm{0}
\end{equation}
where $\bm{u} \in \mathbb{R}^{h}$ is the primary vector field of interest, $\bm{x} \in \Omega \subseteq \mathbb{R}^{dim}$ is the spatial field with $\Omega$ being the domain of definition of the PDE, $t \in \mathbb{R}$ is the time coordinate, $\tilde{\bm{G}}$ is an operator that maps into $\mathbb{R}^{N_{G}}$ and $\bm{\mathcal{F}}$ is a set of fuzzy fields $\mathcal{F} = \lbrace (F_{1}(\bm{x},t), \mu_{1}) \ldots, (F_{s}(\bm{x},t), \mu_{s})  \rbrace$ all defined by a certain convex membership function.
At a given $\alpha$-cut value we are able to transform this fuzzy PDE into an interval PDE defined as
\begin{equation}\label{eq::IPDE}
    \bm{G}(\bm{x},t, \bm{u}, \mathcal{P}_{\alpha}) = \bm{0}
\end{equation}
where $\bm{G}$ still maps into $\mathbb{R}^{N_{G}}$ and $\mathcal{P}_{\alpha}$ is now the set of $\alpha$-cut interval fields found at $\mu(F_{i}(\bm{x},t))=\alpha$, $\forall i=1,\ldots,s$
\begin{equation}
    \mathcal{P}_{\alpha} = \lbrace P_{1}(\bm{x},t), \ldots, P_{s}(\bm{x},t)  \rbrace
\end{equation}
and each $P_{i}(\bm{x},t)$ is defined by $P_{i}(\bm{x},t) \in P_{i}^{I}(\bm{x},t) = \left[P_{i}^{L}(\bm{x},t), P_{i}^{U}(\bm{x},t) \right]$.
The solution $\tilde{\bm{u}}(\bm{x}, t,\mathcal{P}_{\alpha})$ to problem (\ref{eq::IPDE}) is complicated and can typically only be obtained approximately.
Consider the time and spatial domain to be discretized into $n$ and $m$ points respectively. 
In this work we are interested in the two solutions of the problem that are given by $\tilde{\bm{u}}_{min}(\bm{x}, t) = \tilde{\bm{u}}(\bm{x}, t,\mathcal{P}_{min})$ and $\tilde{\bm{u}}_{max}(\bm{x}, t) = \tilde{\bm{u}}(\bm{x}, t,\mathcal{P}_{max})$ where
\begin{equation}
    \begin{aligned}
        \mathcal{P}_{min} &= \min_{\mathcal{P}_{\alpha}^{\star} \in \mathcal{P}_{\alpha}} \sum_{i=1}^{m} \sum_{j=1}^{n} \sum_{k=1}^{h} \tilde{u}_{k}(\bm{x}_{i}, t_{j},\mathcal{P}_{\alpha}^{\star}) \\
        \mathcal{P}_{max} &= \max_{\mathcal{P}_{\alpha}^{\star} \in \mathcal{P}_{\alpha}} \sum_{i=1}^{m} \sum_{j=1}^{n} \sum_{k=1}^{h} \tilde{u}_{k}(\bm{x}_{i}, t_{j},\mathcal{P}_{\alpha}^{\star}) 
    \end{aligned}
\end{equation}
where $\tilde{u}_{k}$ is the $k$-th component of the vector.
Hence, we are interested in the two solutions which correspond to the two cases inside the allowed interval set in which the input parameter values minimize and maximize the sum of the parameter of interest with
\begin{equation}
    \begin{aligned}
    \mathcal{P}_{min}&= \lbrace P_{1}^{min}(\bm{x},t), \ldots, P_{s}^{min}(\bm{x},t) |  P_{i}^{L}(\bm{x},t) \leq P_{i}^{min}(\bm{x},t) \leq P_{i}^{U}(\bm{x},t) , \forall i=1,\ldots, s \rbrace , \\
        \mathcal{P}_{max}&= \lbrace P_{1}^{max}(\bm{x},t), \ldots, P_{s}^{max}(\bm{x},t) |  P_{i}^{L}(\bm{x},t) \leq P_{i}^{max}(\bm{x},t) \leq P_{i}^{U}(\bm{x},t) , \forall i=1,\ldots, s \rbrace .
    \end{aligned}
\end{equation}
In the limit, $\tilde{\bm{u}}_{min}(\bm{x}, t)$ and $\tilde{\bm{u}}_{max}(\bm{x}, t)$ will yield the output interval bounds.
From an engineering point of view these two extreme cases are important because they represent the worst case scenarios of the problem at hand which has significant impact on the design of structures, for example.
\begin{figure}[hb!]
    \centering
    \includegraphics[scale=0.65]{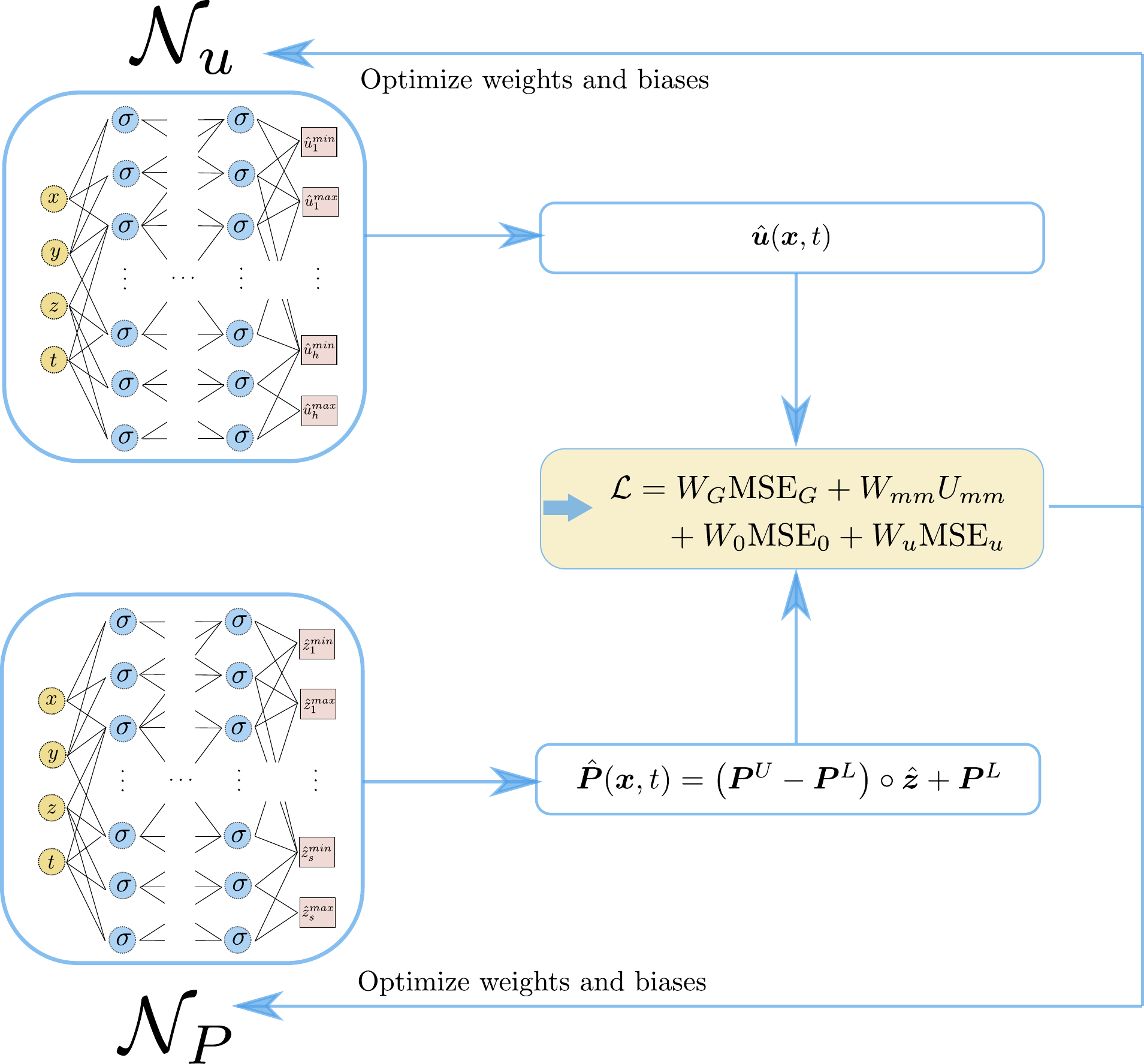}
    \caption{Overview of iPINN framework}
    \label{fig:iPINNFrame}
\end{figure}
\section{Interval and fuzzy physics-informed neural networks}\label{sec::3}
In this section we provide a brief overview of the neural networks formulation.
Neural networks are in general composed of one input, one output and $n_{D}-1$ hidden layers.
Let the weights and biases of the $k^{\text{th}}$ layer be denoted by $\bm{W}^{k} $ and $\bm{b}^{k}$. 
Consider that the $k^{\text{th}}$ hidden layer transfers some output $\bm{x}^{k}$ to the $(k+1)^{\text{th}}$ layers which applies an affine transformation
\begin{equation}\label{eq::NNO}
    \mathcal{L} (\bm{x}^{k})  = \bm{W}^{k+1} \bm{x}^{k} + \bm{b}^{k+1}
\end{equation}
and some activation function $\sigma$ to it.
Since, equation (\ref{eq::NNO}) is applied in every layer of the network an input $\bm{x}$ yields a network output of the form
\begin{equation}
    \hat{\bm{y}}(\bm{x}) = (\mathcal{L}_{k} \circ \sigma \circ \mathcal{L}_{k-1}\circ \cdots \circ \sigma \circ \mathcal{L}_{1})(\bm{x})
\end{equation}
where $\circ$ is a composition operator. 
The goal of neural networks is to find the optimal set of trainable parameters $\bm{\Theta} = \lbrace \bm{W}^{k}, \bm{b}^{k}  \rbrace_{k=1}^{n_{D}}$ such that the network $\mathcal{N}(\bm{\Theta})$ provides the best fit for the input-output mapping.
This is achieved 
by following an optimization procedure defined over some loss function $J(\bm{\Theta} )$
\begin{equation}\label{eq::0_optimization}
    \bm{\Theta}^{\star} = \argmin_{\bm{\Theta}} J(\bm{\Theta} ).
\end{equation}
This problem is typically solved in an iterative manner by employing a stochastic gradient-descent approach. A review on neural network optimization methods is provided in \cite{bottou2018optimization}.
For more information on neural networks we refer to \cite{goodfellow2016deep}. 
The proposed interval physics-informed neural network (iPINN) consists of two separate feedforward neural networks one to approximate the maximum and minimum solution outputs termed $\mathcal{N}_{u}(\bm{\Lambda})$ and one to find the corresponding input fields $\mathcal{N}_{P}(\bm{\Psi})$ where $\bm{\Lambda}$ and $\bm{\Psi}$ are the trainable parameters, see Figure \ref{fig:iPINNFrame}. Both networks define a global shape function over the time-space domain $(\bm{x},t)$ of the problem. We define the final activation function of $\mathcal{N}_{u}(\bm{\Lambda})$ as identity whereas the equivalent function of $\mathcal{N}_{P}(\bm{\Psi})$ is of sigmoidal form. 
This is done to restrict the outputs $\hat{\bm{z}}(\bm{x},t, \bm{\Lambda})$ of $\mathcal{N}_{P}(\bm{\Psi})$ to be between $0$ and $1$ in order to allow for easy scaling of the outputs into the required interval field, i.e.
\begin{equation}
    \begin{aligned}
\hat{\bm{P}}(\bm{x},t) =\left(\bm{P}^{U}(\bm{x},t) - \bm{P}^{L}(\bm{x},t)\right) \circ \hat{\bm{z}}(\bm{x},t, \bm{\Lambda}) + \bm{P}^{L}(\bm{x},t),
    \end{aligned}
\end{equation}
where $\bm{P}^{U}(\bm{x},t) =\begin{bmatrix}
        P_{1}^{U},
        P_{1}^{U},
        \ldots,
        P_{i}^{U},
        P_{i}^{U},
        \ldots,
        P_{s}^{U},
        P_{s}^{U}
        \end{bmatrix}^{T}$ and $\bm{P}^{L}$ is defined equivalently but with the lower limit of the interval.
The output of $\mathcal{N}_{u}(\bm{\Lambda})$  is of the form $\hat{\bm{u}}(\bm{x},t, \bm{\Lambda}) = \begin{bmatrix}
        \hat{u}_{1}^{min},
        \hat{u}_{1}^{max},
        \ldots,
        \hat{u}_{h}^{min},
        \hat{u}_{h}^{max}
        \end{bmatrix}^{T}$.
As typically found in PINN \citep{raissi2019physics}
consider a set of domain training points $\lbrace \bm{X}_{\mathcal{R}}^{i} \rbrace_{i=1}^{N_{\mathcal{R}}}$, boundary training points $\lbrace \bm{X}_{u}^{i} \rbrace_{i=1}^{N_{u}}$ and initial state training points $\lbrace \bm{X}_{0}^{i} \rbrace_{i=1}^{N_{0}}$ and some time discretization $\lbrace t^{i} \rbrace_{i=1}^{N_{t}}$  where $N_{\mathcal{R}}$, $N_{u}$, $N_{0}$ and $N_{t}$ denote the number of points, respectively.
We can then formulate the optimization problem for iPINN in the following way
\begin{equation}\label{eq::WeightsWithLoss}
\begin{aligned}
        (\bm{\Lambda}^{\star}, \bm{\Psi}^{\star}) =  \argmin_{(\bm{\Lambda}, \bm{\Psi})}\quad &W_{G} \, \text{MSE}_{G} \left(\lbrace \lbrace \bm{X}_{\mathcal{R}}^{i} \rbrace_{i=1}^{N_{\mathcal{R}}}, \lbrace t^{j} \rbrace_{j=1}^{N_{t}} \rbrace,  (\bm{\Lambda}, \bm{\Psi}) \right) +  W_{mm} \, U_{mm}\left(\lbrace \lbrace \bm{X}_{\mathcal{R}}^{i} \rbrace_{i=1}^{N_{\mathcal{R}}}, \lbrace t^{j} \rbrace_{j=1}^{N_{t}} \rbrace, \bm{\Lambda}  \right) \\&+ W_{0} \, \text{MSE}_{0} \left(\lbrace \lbrace \bm{X}_{t}^{i} \rbrace_{i=1}^{N_{t}},   t^{0}   \rbrace,\bm{\Lambda} \right) +W_{u} \, \text{MSE}_{u} \left(\lbrace \lbrace \bm{X}_{u}^{i} \rbrace_{i=1}^{N_{u}}, \lbrace t^{j} \rbrace_{j=1}^{N_{t}} \rbrace, \bm{\Lambda} \right)
\end{aligned}
\end{equation}
where 
\begin{equation}\label{eq::MainLoss}
    \begin{aligned}
       \text{MSE}_{G} \left(\lbrace \lbrace \bm{X}_{\mathcal{R}}^{i} \rbrace_{i=1}^{N_{\mathcal{R}}}, \lbrace t^{j} \rbrace_{j=1}^{N_{t}} \rbrace,  (\bm{\Lambda}, \bm{\Psi}) \right) &= \frac{1}{N_{\mathcal{R}} N_{t}} \sum_{i=1}^{N_{\mathcal{R}} }\sum_{j=1}^{N_{t} } \sum_{k=1}^{N_{G}} \abs{ G_{k}(\bm{X}_{\mathcal{R}}^{i},t^{j}, \hat{\bm{u}}(\bm{\Lambda}), \hat{\bm{P}}(\bm{\Psi}))}^{2}\\
       U_{mm}\left(\lbrace \lbrace \bm{X}_{\mathcal{R}}^{i} \rbrace_{i=1}^{N_{\mathcal{R}}}, \lbrace t^{j} \rbrace_{j=1}^{N_{t}} \rbrace, \bm{\Lambda}  \right)&= \sum_{i=1}^{N_{\mathcal{R}} }\sum_{j=1}^{N_{t} } \sum_{k=1}^{h} \hat{u}_{k}^{min}(\bm{X}_{R}^{i} , t^{j}) 
       - \sum_{i=1}^{N_{\mathcal{R}} }\sum_{j=1}^{N_{t} } \sum_{k=1}^{h} \hat{u}_{k}^{max}(\bm{X}_{R}^{i} , t^{j}) \\
       \text{MSE}_{0} \left(\lbrace \lbrace \bm{X}_{0}^{i} \rbrace_{i=1}^{N_{0}},   t^{0}   \rbrace,\bm{\Lambda} \right) &= \frac{1}{N_{0}} \sum_{i=1}^{N_{0} } \sum_{k=1}^{h} \abs{\hat{u}_{k}^{min}(\bm{X}_{0}^{i} , t^{0}) - u_{k}(\bm{X}_{0}^{i} , t^{0}) }^{2} \\&  + \abs{\hat{u}_{k}^{max}(\bm{X}_{0}^{i} , t^{0})- u_{k}(\bm{X}_{0}^{i} , t^{0}) }^{2} \\
       \text{MSE}_{u} \left(\lbrace \lbrace \bm{X}_{u}^{i} \rbrace_{i=1}^{N_{u}}, \lbrace t^{j} \rbrace_{j=1}^{N_{t}} \rbrace, \bm{\Lambda} \right) &= \frac{1}{N_{u} N_{t}} \sum_{i=1}^{N_{u} } \sum_{j=1}^{N_{t} }\sum_{k=1}^{h} \abs{\hat{u}_{k}^{min}(\bm{X}_{u}^{i} , t^{j}) - u_{k}(\bm{X}_{u}^{i} , t^{j})}^{2} \\&    + \abs{\hat{u}_{k}^{max}(\bm{X}_{u}^{i} , t^{j})- u_{k}(\bm{X}_{u}^{i} , t^{j})}^{2}
    \end{aligned}
\end{equation}
and $W_{G}$, $W_{mm}$, $W_{0}$ and $W_{u}$ are user-chosen weights.

To better clarify the intuition behind this particular choice for the loss function, an explanation of its constituents is provided next. $\text{MSE}_{G}$ is the term that needs to be minimized in order to ensure that the solution $\bm{u}$ accurately satisfies the PDE at the interior of its domain $\Omega$. In a similar fashion, the purpose of $\text{MSE}_{0}$ and $\text{MSE}_{u}$ is to impose the initial and boundary conditions, respectively, on the solution. Lastly, $U_{mm}$ is the additional term in the iPINN framework, whose minimization eventually yields the extrema $\hat{\bm{u}}^{min}$ and $\hat{\bm{u}}^{max}$. 

\begin{figure}[ht!]
    \centering
    \includegraphics[scale=0.8]{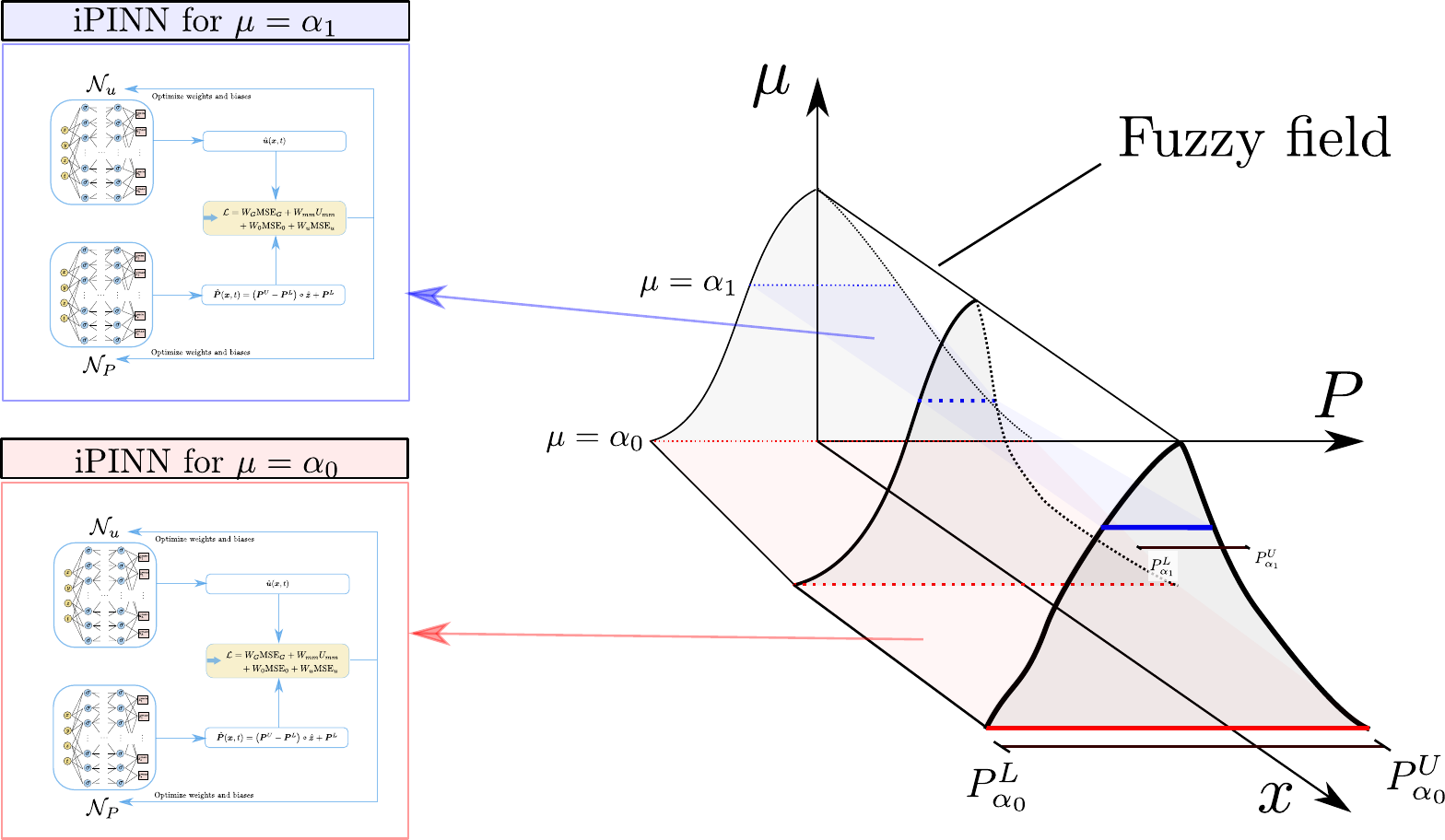}
    \caption{Overview of fPINN framework as a collection of iPINNs}
    \label{fig:fPINNFrame}
\end{figure}

A fuzzy physics-informed neural network (fPINN) can then be understood as a collection of iPINNs, see Figure \ref{fig:fPINNFrame}. 
After minimization of the loss function of eq. (\ref{eq::MainLoss}), the trainable parameters of $\mathcal{N}_{u}$ and $\mathcal{N}_{P}$ are obtained. The outputs of these networks yield approximate bounds to iPDE solutions ($\mathcal{N}_{u}$) as well as the corresponding input field values that lead to these extreme values ($\mathcal{N}_{P}$). It needs to be highlighted that in contrast to FEM no correlation length specification of the input fields needs to be assumed and Monte-Carlo-like simulations are also not required.
Additionally, we retain all major advantages of the PINN approach, i.e. the problem is solved meshfree and inverse problems are straightforwardly set up and solvable. For more information on PINNs we refer to \cite{raissi2019physics}. One slightly negative aspect of this approach is that we employ two network outputs for $\mathcal{N}_{u}$ and $\mathcal{N}_{P}$ for each input field which for a larger number of input fields might be restrictive. However this limitation is dependent on the problem complexity and its dimensionality. The presented approach is somewhat related to the stochastic PINN framework proposed in \cite{zhang2019quantifying}. However their main concern is probabilistic uncertainty quantification which is not the focus of this paper.%

\section{Application}\label{sec::4}
In the following, the proposed formulation is tested and studied on a set of problems of increasing complexity. First, a simple educational problem is considered, involving an interval parameter and fuzzy parameter, and thereafter, we attempt to find the solution to a more complex one-dimensional structural problem involving two simultaneously applied interval fields. Lastly, we study the performance of the framework on a time-dependent interval PDE problem.
The fPINN/iPINN formulation was implemented in
Pytorch \citep{NEURIPS2019_9015} and the network parameters were optimized using the
Adam optimizer \citep{kingma2014adam}\footnote{Codes will be made public under \url{https://github.com/FuhgJan/intervalAndFuzzyPINN} after acceptance of this paper.}. The FEM results for the second problem were obtained using the FEniCS framework \citep{AlnaesBlechta2015a}.
\subsection{Introductory problem: Function of an interval and fuzzy parameter}
\begin{figure}[ht!]
\centering
\includegraphics[scale=0.3]{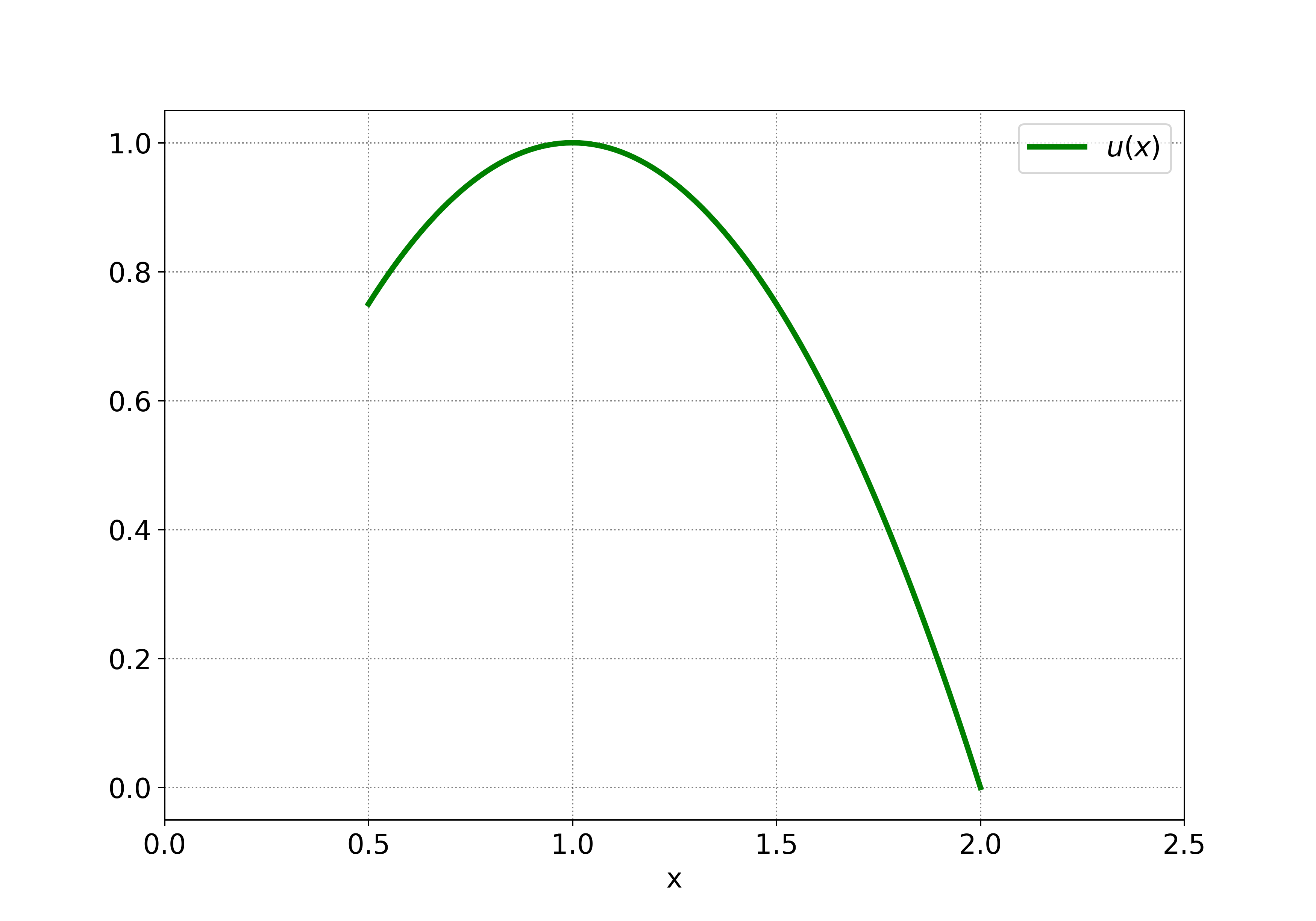} 
\caption{Dependence of primary variable over elements of interval variable}\label{fig:ToyPrimary}
\end{figure}

Consider the following simple functional problem with an interval parameter
\begin{equation}
    u^{I}(x^{I}) = x^{I} (2 - x^{I}), 
\end{equation}
where the interval parameter is constrained by
\begin{equation}
    x^{I} = \left[ 0.5 ,2.0 \right].
\end{equation}
The primary output variable $u^{I}(x)$ is plotted over the elements of the interval input in Figure \ref{fig:ToyPrimary}. It can be seen that $u$ is non-monotonically dependent on $x$. As a result the output bounds $\hat{u}^{min}$ and $\hat{u}^{max}$ are not just resulting from endpoint combinations, i.e. just the evaluation of $x = x^{L} =  0.5$ and $x = x^{U} = 2.0$ is not enough to obtain the output bounds for more information see e.g. \citep{ZADEH1965338}. In fact as seen in Figure \ref{fig:ToyPrimary} the bounds of the primary variable $\hat{u}^{min}= 0.0$ and $\hat{u}^{max}= 1.0$ are results of the input interval values $x^{min}=2.0$ and $x^{max}=1.0$ respectively.

\begin{figure}[ht]
\begin{subfigure}[b]{0.5\linewidth}
\centering
\includegraphics[scale=0.27]{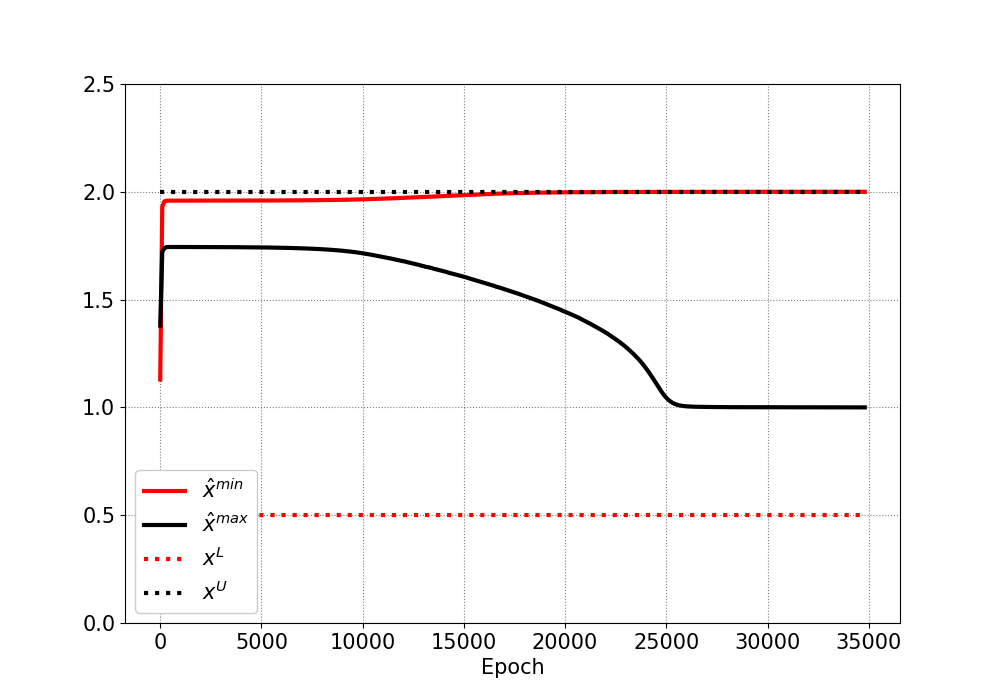} 
\caption{Predicted and constraining input field}\label{fig:toyExampleInput}
\end{subfigure}%
\begin{subfigure}[b]{.5\linewidth}
\centering
\includegraphics[scale=0.27]{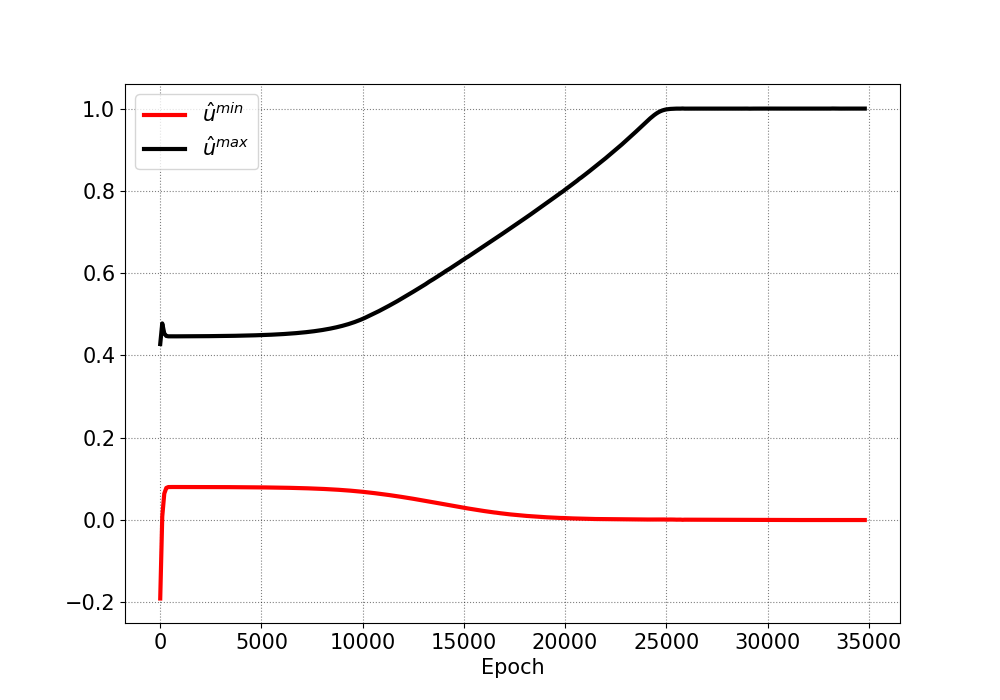} 
\caption{Predicted primary output bounds}\label{fig:toyExampleOutput}
\end{subfigure}
\caption{Results of introductory example. Outputs of both trained networks over training iterations.}\label{fig:}
\end{figure}

To highlight how iPINN is able to solve this problem we consider both $\mathcal{N}_{u}$ and $\mathcal{N}_{P}$ consisting of 2 hidden layers with 20 neurons. The activation functions are chosen to be the hyperbolic tangent function. These choices are arbitrary and are not the result of any hyperparameter optimization, i.e. they do not represent any special network setup to the best of the author's knowledge. Both networks have $1$ input and $2$ outputs.
When dealing with interval fields and iPDEs the spatial and temporal values are the inputs. Since in this simple problem the primary variable is independent of space and time, a constant arbitrary value can be chosen as the input which has no effect on the solution (we chose $1.0$).
The training is conducted with a learning rate of $1e-3$ and is stopped after $35,000$ iterations. High emphasis is put on the fulfilment of the residual $u^{I}(x^{I}) - x^{I} (2 - x^{I}) = 0$. Hence the loss weights in eq. (\ref{eq::WeightsWithLoss}) are chosen to be $W_{G}=100,000$ and $W_{mm}=1$. Since there are no initial and boundary conditions we do not need to account for those terms in the loss function. \\
Figure \ref{fig:toyExampleInput} shows the constraining limit values of the input interval as well as the predicted input values that lead to the output interval bounds, see Figure \ref{fig:toyExampleOutput}, over the training process. It can be seen that both the input values, as well as the output bounds are accurately predicted after around 25,000 epochs.

\begin{figure}[ht!]
\centering
\includegraphics[scale=1.0]{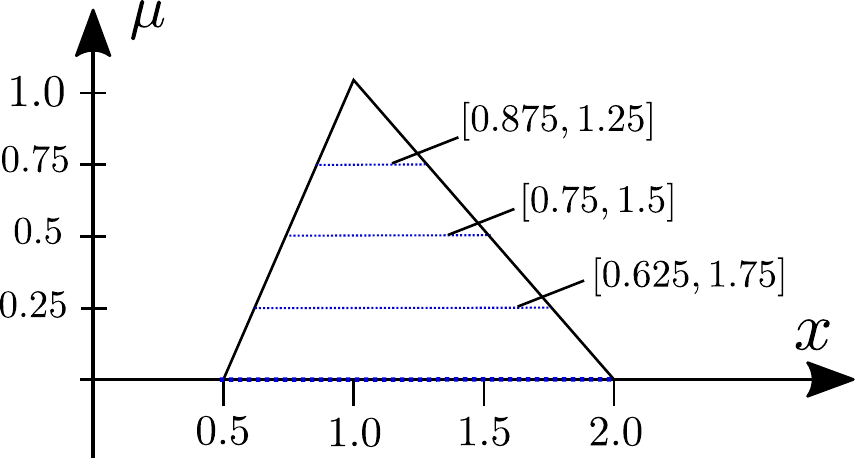} 
\caption{Triangular fuzzy input with interval values at $\alpha$-cuts}\label{fig:ToyFuzzyInp}
\end{figure}

Instead of an interval number, consider the function to be dependent on a fuzzy number
\begin{equation}
    u^{F}(x^{F}) = x^{F} (2 - x^{F}). 
\end{equation}
Let $x^{F}$ be a triangular fuzzy number of the form given in Figure \ref{fig:ToyFuzzyInp} where $\mu(x=1)= 1.0$. We can see that the previously investigated problem coincides with the $\alpha$-cut at $\mu=0$. In order to obtain an approximate description of the output fuzzy bounds we consider 4 more $\alpha$-cuts at $0.25, 0.5, 0.75$ and $1.0$. The interval bounds for these cuts are also shown in Figure \ref{fig:ToyFuzzyInp}.

\begin{figure}
\begin{subfigure}[b]{0.5\linewidth}
\centering
\includegraphics[scale=0.27]{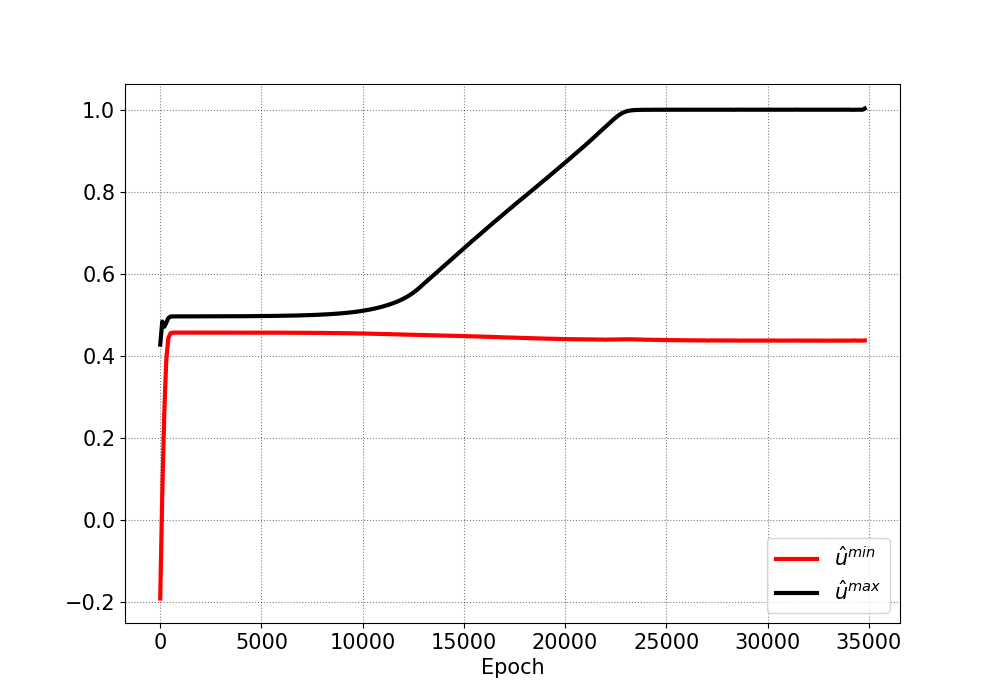} 
\caption{$\alpha$- cut at $\mu=0.25$}
\end{subfigure}%
\begin{subfigure}[b]{.5\linewidth}
\centering
\includegraphics[scale=0.27]{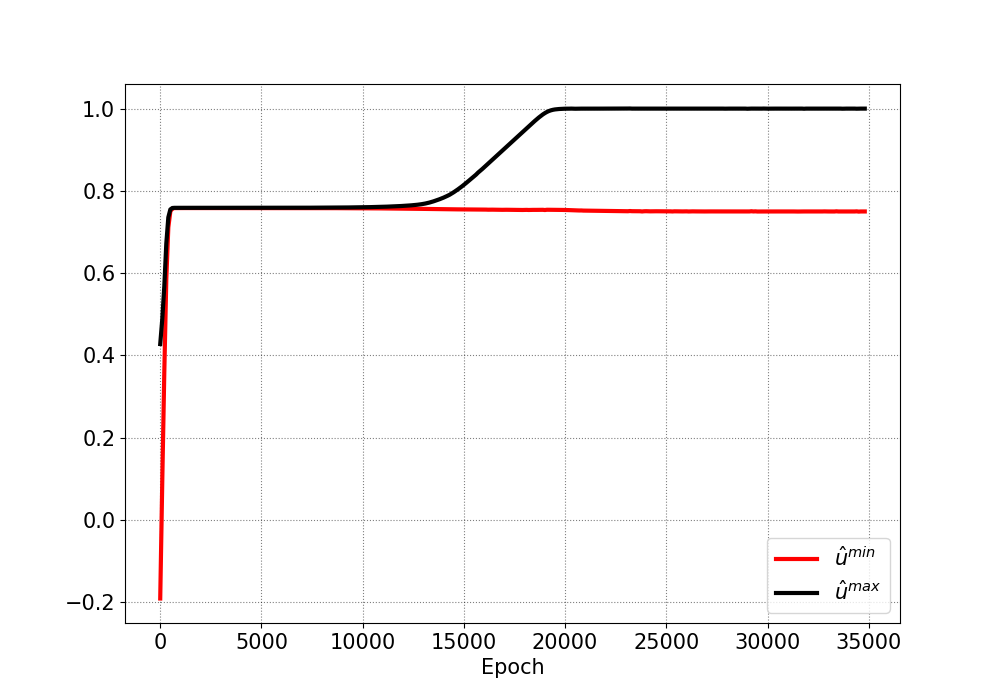} 
\caption{$\alpha$- cut at $\mu=0.5$}
\end{subfigure}
\begin{subfigure}[b]{0.5\linewidth}
\centering
\includegraphics[scale=0.27]{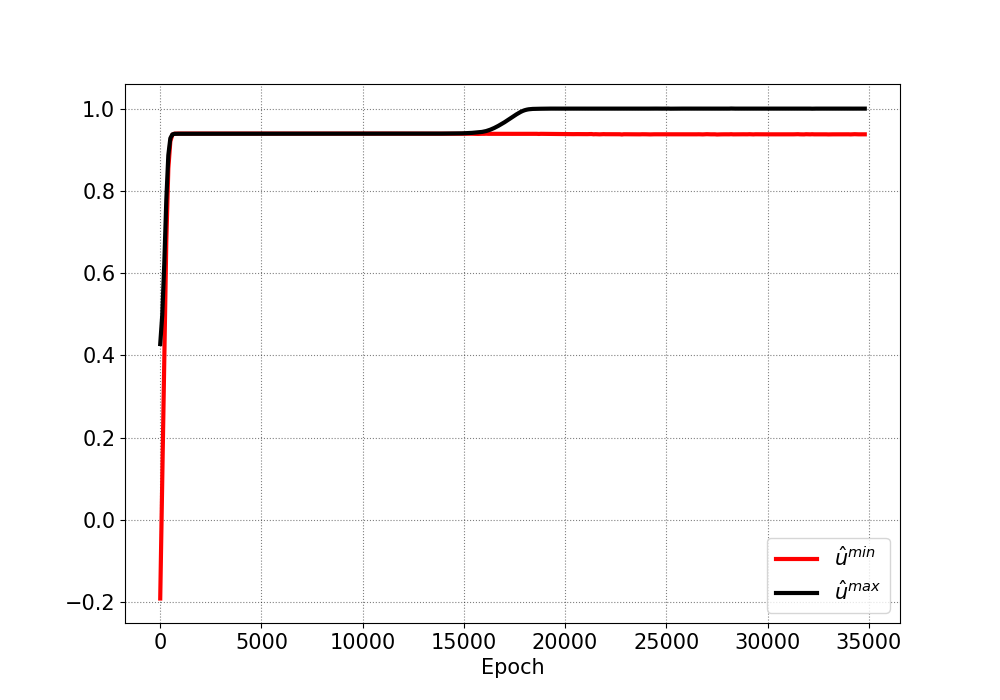} 
\caption{$\alpha$- cut at $\mu=0.75$}
\end{subfigure}%
\begin{subfigure}[b]{.5\linewidth}
\centering
\includegraphics[scale=0.27]{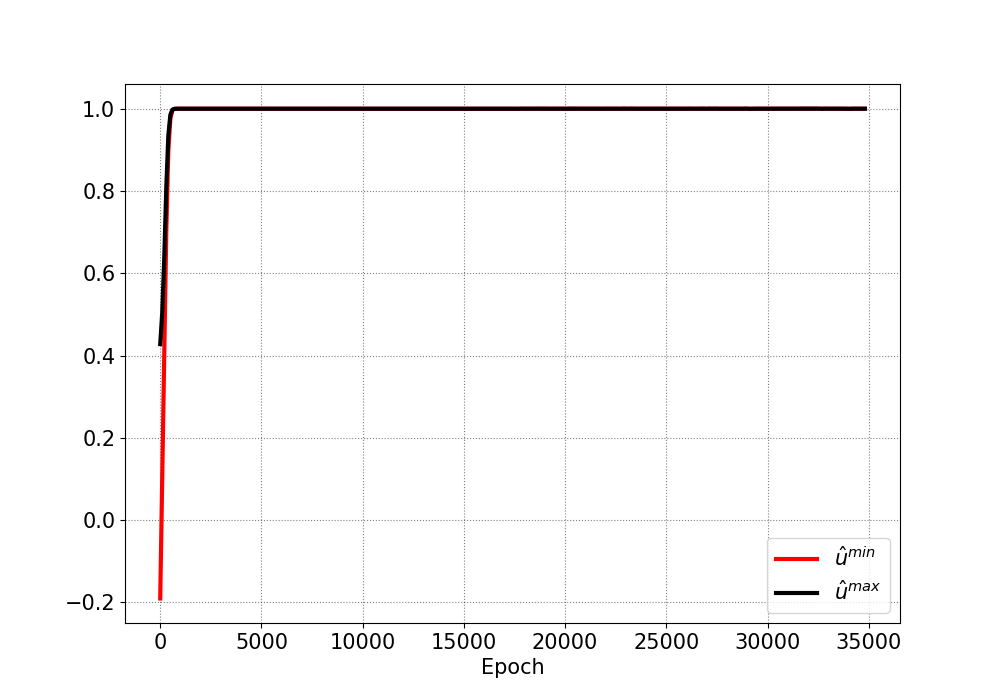} 
\caption{$\alpha$- cut at $\mu=1.0$}
\end{subfigure}
\caption{fPINN predicted bounding outputs for different $\alpha$-cuts over training process}\label{fig:ToyFuzzyOutTrain}
\end{figure}

Following the fPINN concept, we build an iPINN for each new $\alpha$-cut value, i.e. resulting in 4 new iPINNs. Only changing the considered bounding interval in comparison to the previous problem and leaving the rest of the networks unchanged (hyperparameter, learning rate, etc.) the predicted output bounding intervals of the subproblems over the training process are shown in Figure \ref{fig:ToyFuzzyOutTrain}.

\begin{figure}[ht!]
\centering
\includegraphics[scale=1.0]{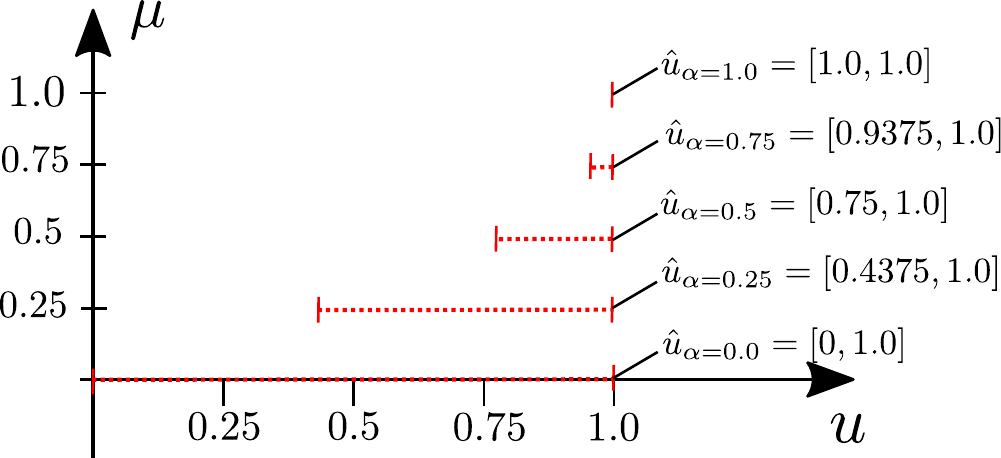} 
\caption{Fuzzy (interval) output values for considered $\alpha$-cuts obtained using fPINN framework}\label{fig:ToyFuzzyOut}
\end{figure}
Hence, using these results for the considered $\alpha$-cuts we can make an approximation of the output fuzzy number by knowing the intervals at each $\alpha$-cut. Therefore, using fPINN we are able to obtain the fuzzy number estimations as shown in Figure \ref{fig:ToyFuzzyOut} which exactly correspond to the analytical solutions that can be obtained trivially. An interesting side note is that the output fuzzy number seen in Figure \ref{fig:ToyFuzzyOut} appears to be no longer triangular.

\clearpage
\subsection{One-dimensional bar with two interval fields}
\begin{figure}
    \centering
    \includegraphics[scale=1.0]{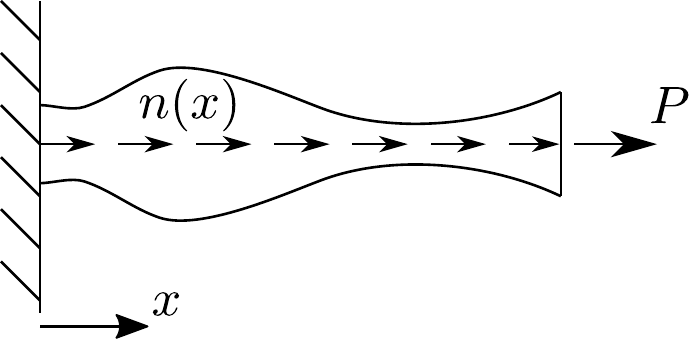}
    \caption{1D bar with area and Young's modulus as interval fields}
    \label{fig:1dBar}
\end{figure}
\begin{figure}
\begin{subfigure}[b]{0.5\linewidth}
\centering
\includegraphics[scale=0.27]{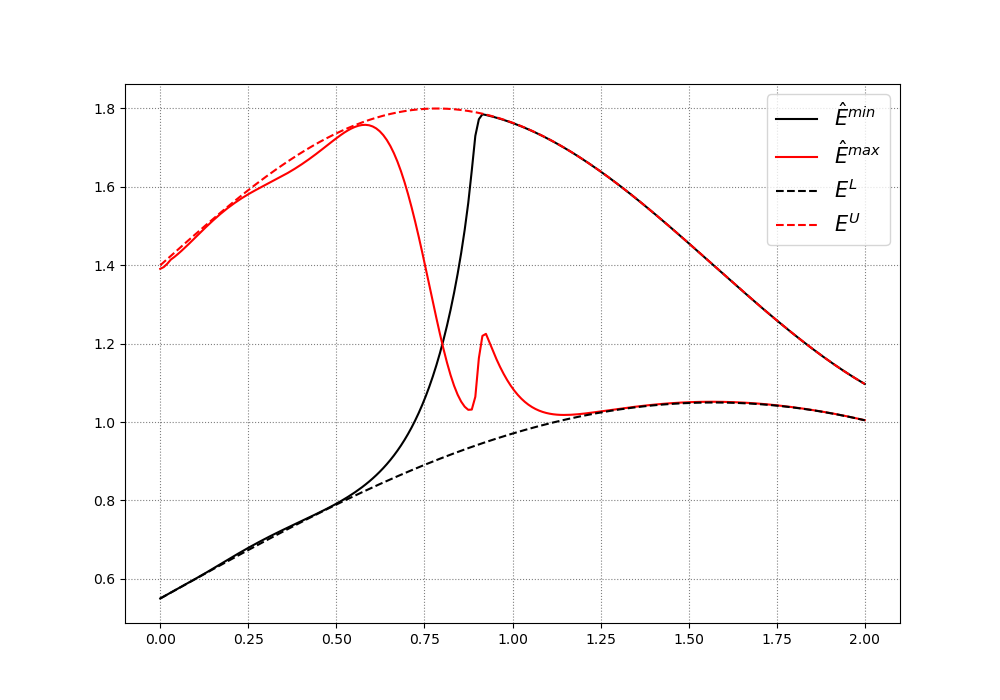} 
\caption{}
\end{subfigure}%
\begin{subfigure}[b]{.5\linewidth}
\centering
\includegraphics[scale=0.27]{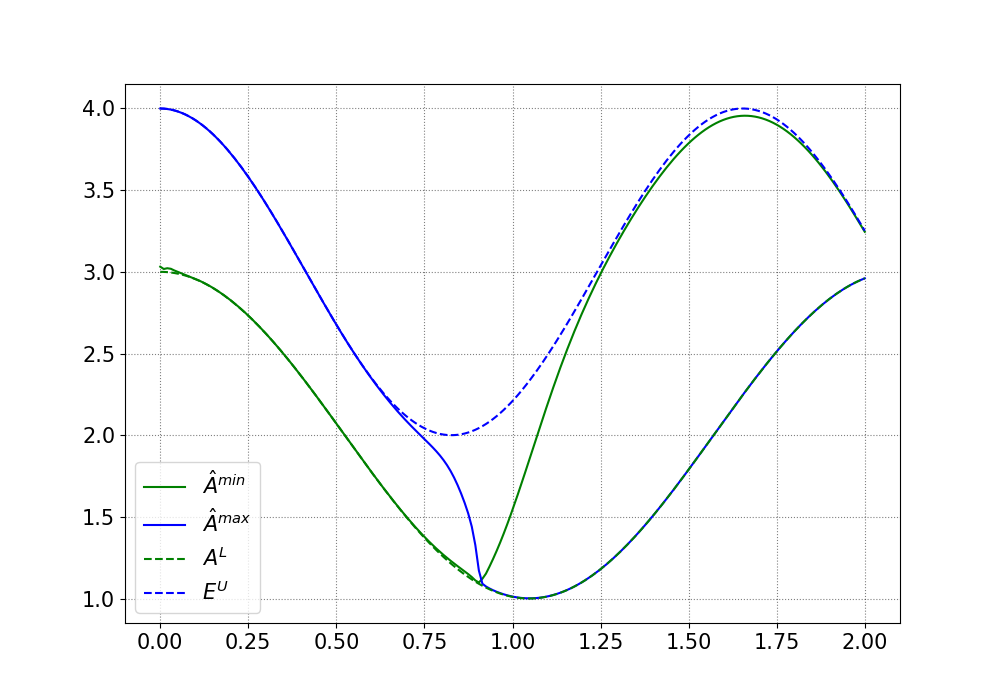} 
\caption{}
\end{subfigure}
\caption{Predicted and constraining input fields, (a) Young's modulus field and (b) area field}\label{fig:InputFields}
\end{figure}
\begin{figure}
\begin{subfigure}[b]{0.5\linewidth}
\centering
\includegraphics[scale=0.27]{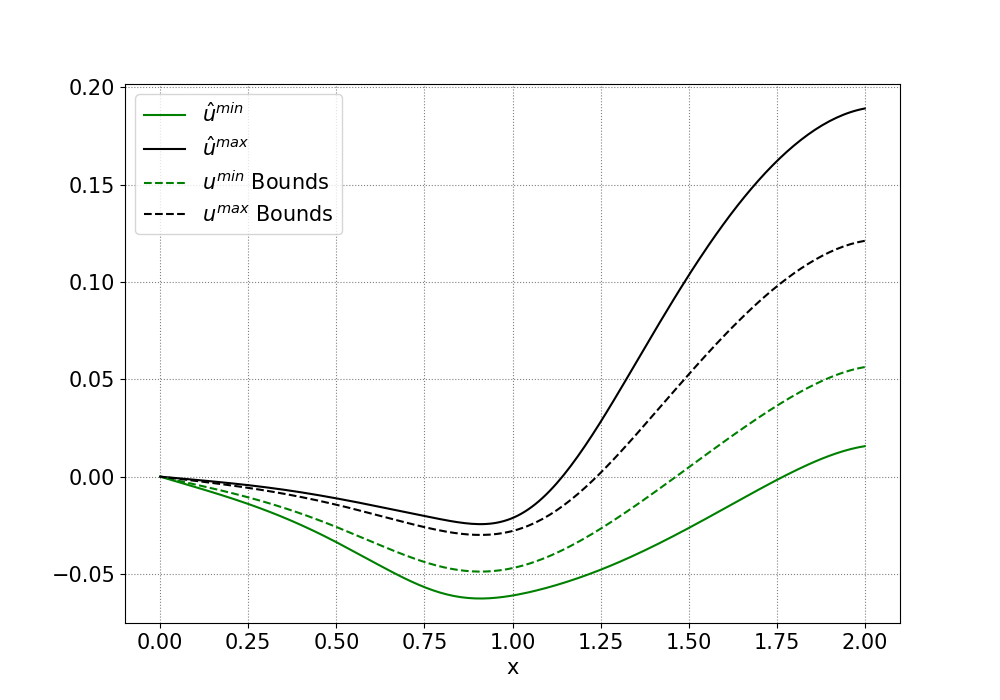}
\caption{}\label{fig::Primaryand4}
\end{subfigure}%
\begin{subfigure}[b]{.5\linewidth}
\centering
\includegraphics[scale=0.27]{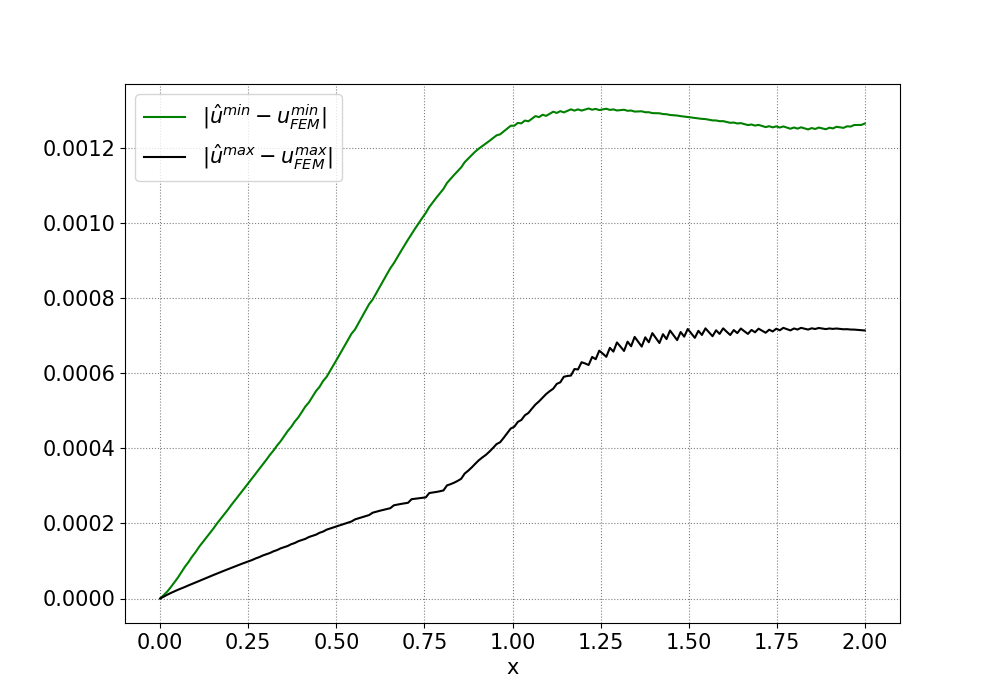}
\caption{}\label{fig::PrimaryError}
\end{subfigure}
\caption{Comparisons and output of primary output. (a) $\mathcal{N}_{u}$ output and maximum and minimum FEM outputs\\ using 4 combination of interval limits, and (b) Error between FEM and the $\mathcal{N}_{u}$ output for $\mathcal{N}_{P}$\\input field}\label{fig:PrimaryFields}
\end{figure}
\begin{figure}[b!]
\begin{subfigure}[b]{0.5\linewidth}
\centering
\includegraphics[scale=0.27]{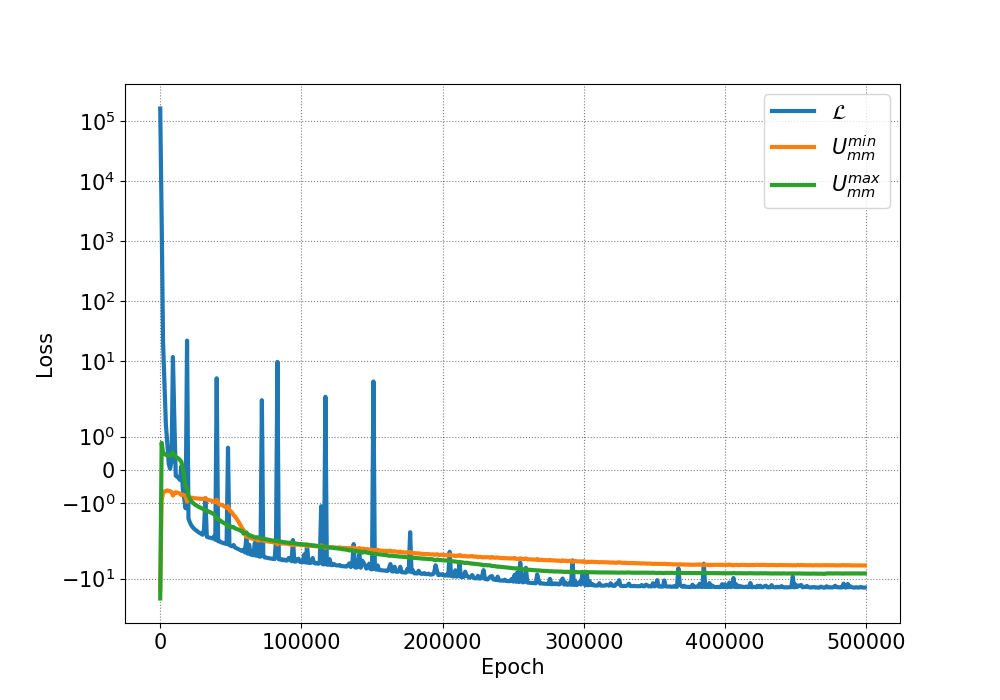} 
\caption{}
\end{subfigure}%
\begin{subfigure}[b]{.5\linewidth}
\centering
\includegraphics[scale=0.27]{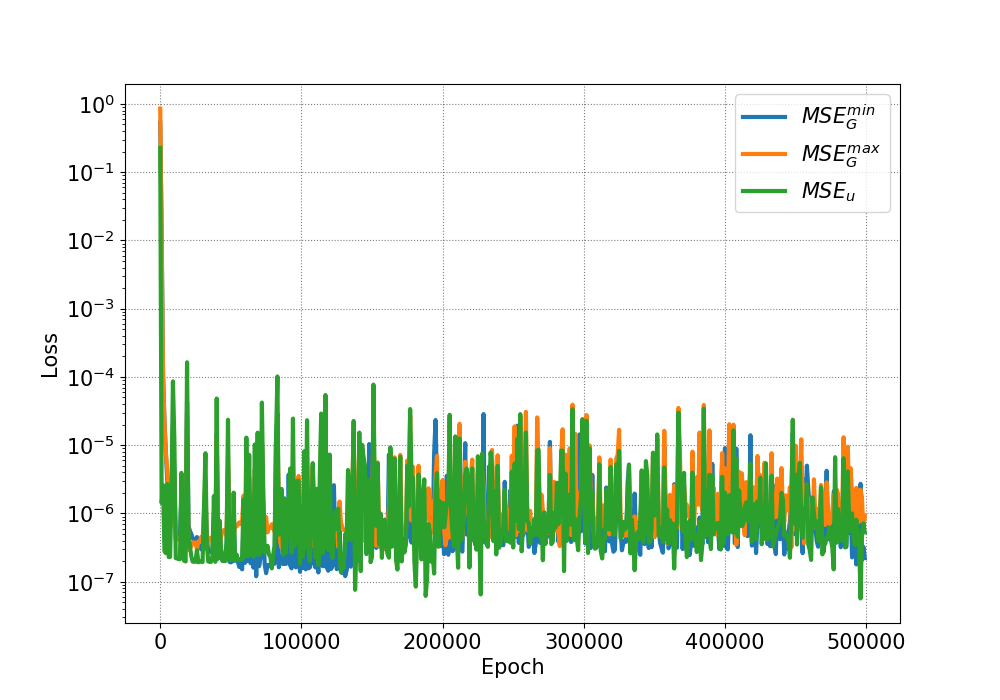} 
\caption{}
\end{subfigure}
\caption{Losses over the training process. (a) Non-residual losses, and (b) Residual losses}\label{fig:Losses}
\end{figure}

In order to test the capability of the proposed approach when dealing with a more realistic structural problem with interval fields, we study the iPINN on a 1D bar problem, see Figure \ref{fig:1dBar}.
The relevant interval ordinary differential equation is given by
\begin{equation}
    \frac{\partial}{\partial x} \left( E^{I}(x) A^{I}(x) \frac{\partial u}{\partial x}\right) + n(x) = 0
\end{equation}
where $n(x)=\cos(3 x) x$ in $N/m$, $E^{I}$ in $N/m^{2}$ and $A^{I}$ in $m^{2}$ are the spatially dependent interval fields of the Young's modulus and the area which are defined as
\begin{equation}
    E^{I}(x) \begin{cases}
    \geq E^{L}(x) &= 0.5 \sin(x) + 0.55  \\
    \leq E^{U}(x) &= 0.4 \sin(2 x) + 1.4
    \end{cases} 
\end{equation}
and 
\begin{equation}
    A^{I}(x) \begin{cases}
    \geq A^{L}(x) &= \cos(3 x) + 2.0\\
    \leq A^{U}(x) &= \cos(3.8 x) + 3.0.
    \end{cases} 
\end{equation}
The boundary conditions are given as 
\begin{equation}\label{eq::Bound}
    u(x=0) = 0, \qquad   \frac{\partial u}{\partial x}\bigg\rvert_{x=L} = \frac{P}{E(x=L) A(x=L)}
\end{equation}
where $P=0.1N$ and where the length of the bar is $L=2m$.
The training is conducted with a learning rate of $1e-4$. The network of the primary solution network $\mathcal{N}_{u}$ consists of 4 hidden layers with 40 neurons each while the input interval network $\mathcal{N}_{P}$ has 5 hidden layers a 50 neurons. The authors can not find anything special about these hyperparameter values. They are employed without any immediate considerations.
We use 200 equidistant points inside the one-dimensional training domain. The FEM mesh consists of 200 elements.
During the studies of this work it was found that high-emphasis needs to be placed on the fulfillment of the residual, i.e. $W_{G}=100,000$ in in eq. (\ref{eq::WeightsWithLoss}) whereas the other weights are less influential for this study $W_{mm}=W_{0}=W_{u}=1$. The network is trained for 500,000 epochs.
Due to the unique formulation of iPINN we are able to obtain the input fields leading to the bounds of the interval displacement field as a simple byproduct of the training approach, see Figure \ref{fig:InputFields}.
It can be seen that the predicted input fields are effectively equal to the constraining interval values for a majority of the computational domain. However both input field predictions "switch" from being closely aligned to the lower bound to being closely aligned to the upper bound value and vice versa. This is an interesting observation which is not trivially predictable beforehand.\\
The bounds of the interval displacement field as well as the bounds of the FEM solutions using the four possible combinations of input interval limits ($E^{L}-A^{L}, E^{U}-A^{L}, E^{L}-A^{U}, E^{U}-A^{U}$) are shown in Figure \ref{fig::Primaryand4}. The absolute error between the FEM solution and the iPINN output for the input fields obtained through $\mathcal{N}_{P}$ are shown in Figure \ref{fig::PrimaryError} which proves that the approach is able to effectively solve the problem in a similar accuracy to FEM. It can be seen that the "switch" of the input fields happens around the global minimum of the primary variable.
The total loss and all individual losses over the training process are shown in Figure \ref{fig:Losses}. It can be seen that the losses decline in a satisfactory manner and that especially the residual losses ($MSE_{G}$ and the boundary loss) are fulfilled quite accurately due to the choice of the loss weights as described above.
In contrast to FEM approaches we do not need to sample from the input fields, i.e. we do not need to assume any spatial correlations of the interval fields. Furthermore, we also do not need to conduct  Monte-Carlo sampling as a means to obtain approximations for the output bounds.

\clearpage
\subsection{Interval nonlinear partial differential equation}
Finally, we investigate the performance of the iPINN framework on the following time-dependent PDE
\begin{equation}\label{eq::nonPDE}
    \frac{\partial u}{\partial t} = 0.01 u \frac{\partial^{2} u}{\partial x^{2}} - k^{I}(x,t) u^{3} + \left( k^{I}(x,t) \right)^{3}
\end{equation}
with the boundary conditions
\begin{equation}
    \begin{aligned}
        u(-1,t) &= u(1,t) = 0 \\
        u(x,0) &= 1-x^{2}.
    \end{aligned}
\end{equation}
Here, $k^{I}(x,t)$ represents a temporally and spatially varying interval field. 
Assume the spatial components to be defined by $x \in [-1,1]$ and the time to be restricted to $t \in [0,1]$.

\begin{figure}
\begin{subfigure}{0.5\linewidth}
    \includegraphics[scale=0.35]{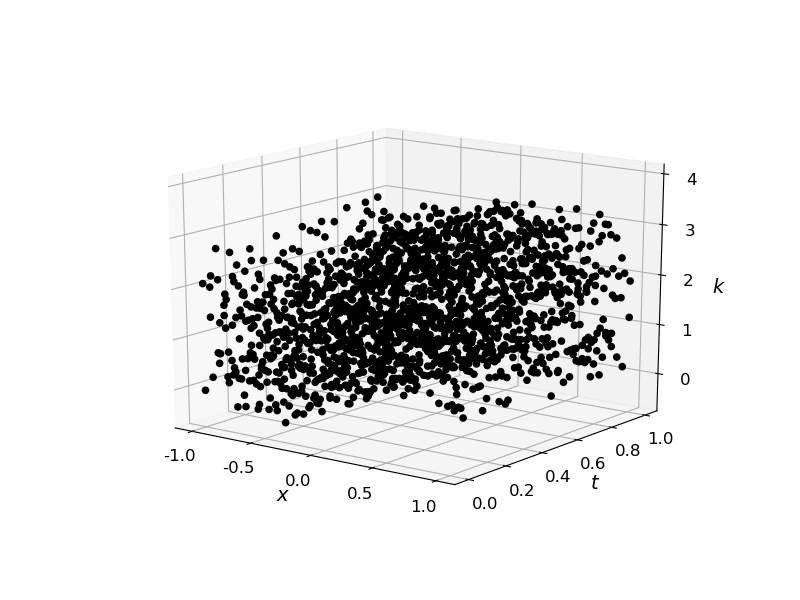}
    \caption{}\label{fig::experimentalPoints}
\end{subfigure}
\begin{subfigure}{0.5\linewidth}
    \includegraphics[scale=0.35]{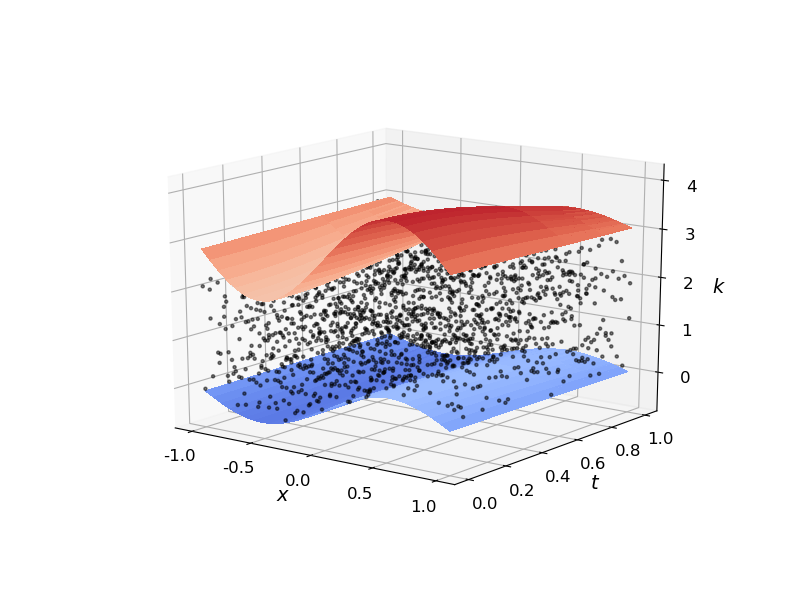}
    \caption{}\label{fig::experimentalPointsFitted}
\end{subfigure}
    \caption{Nonlinear interval PDE problem. (a) "Experimentally" obtained points and (b) fitted lower and upper limits.}
    \label{fig:my_label}
\end{figure}

Suppose we are able to measure $k(x,t)$ the temporal and spatial variation of the field of interest. These experimentally obtained points are shown in Figure \ref{fig::experimentalPoints}. Furthermore, assume that from these points we are able to define lower and upper bounding limits of $k(x,t)$ (Figure \ref{fig::experimentalPointsFitted}) given by
\begin{equation}
     k^{I}(x,t)  \begin{cases}   \geq k^{L}(x,t) &= 0.5 \sin(3 x) \cos(t),\\
    \leq k^{U}(x,t) &= \sin(3x) \cos^{2}(t) + 3.
    \end{cases}
\end{equation}

Both, the input interval network $\mathcal{N}_{P}$ and the primary solution network $\mathcal{N}_{u}$ are chosen to have 3 hidden layers with 40 neurons.
Here, we iterate again that we can not find anything special about these hyperparameter values. They are employed without any immediate considerations.
A learning rate of $1e-4$ for the ADAM optimizer is employed.
We discretize the time domain with $50$ and the spatial domain with $125$ equidistant points.
Finally, we choose $W_{G} = W_{0} = 100,000$ and all other weights equal to $1$.

\begin{figure}
\begin{subfigure}{.5\linewidth}
\includegraphics[scale=0.35]{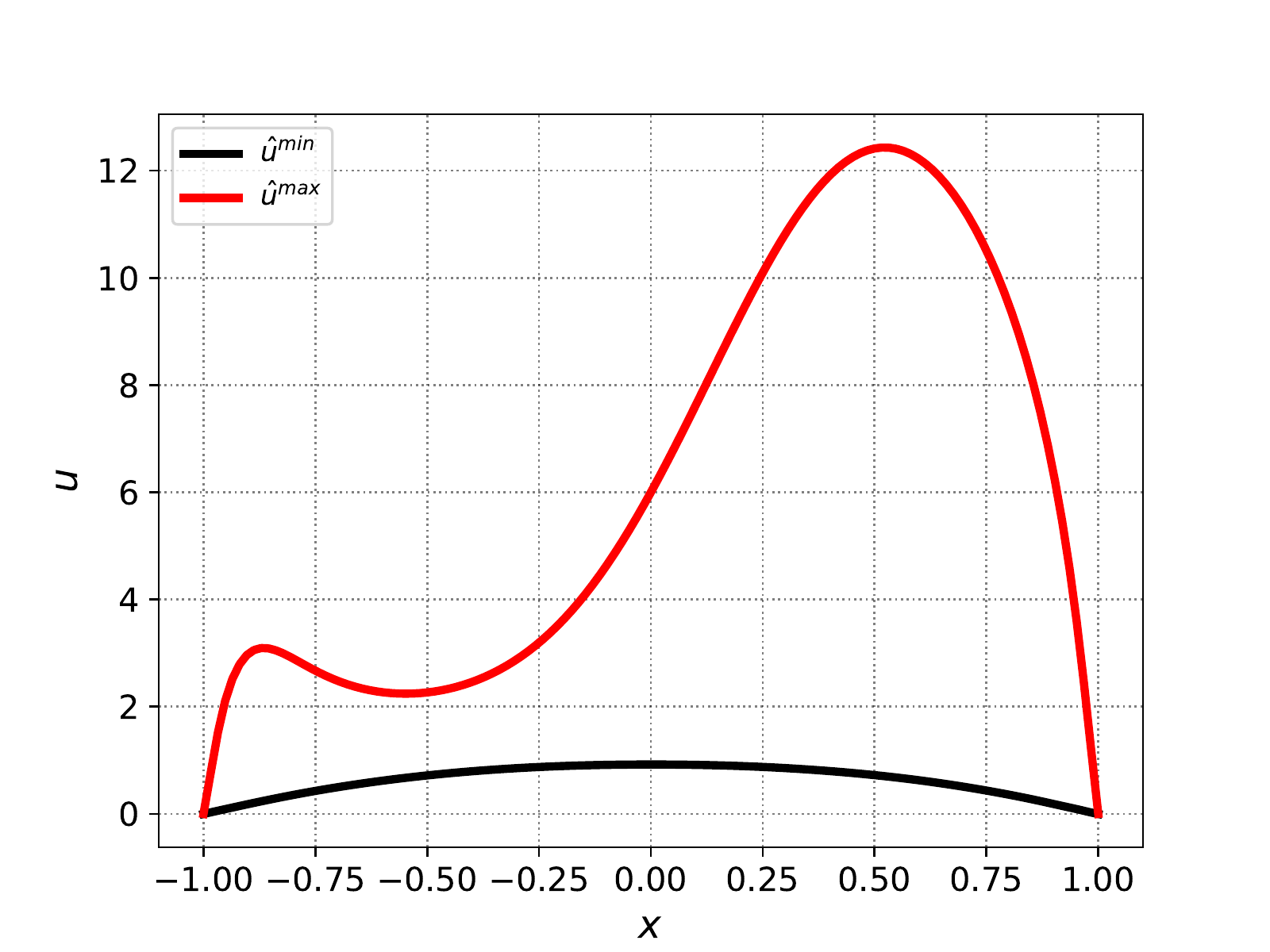}
\caption{$t=0.2s$}
\end{subfigure}
\begin{subfigure}{.5\linewidth}
\includegraphics[scale=0.35]{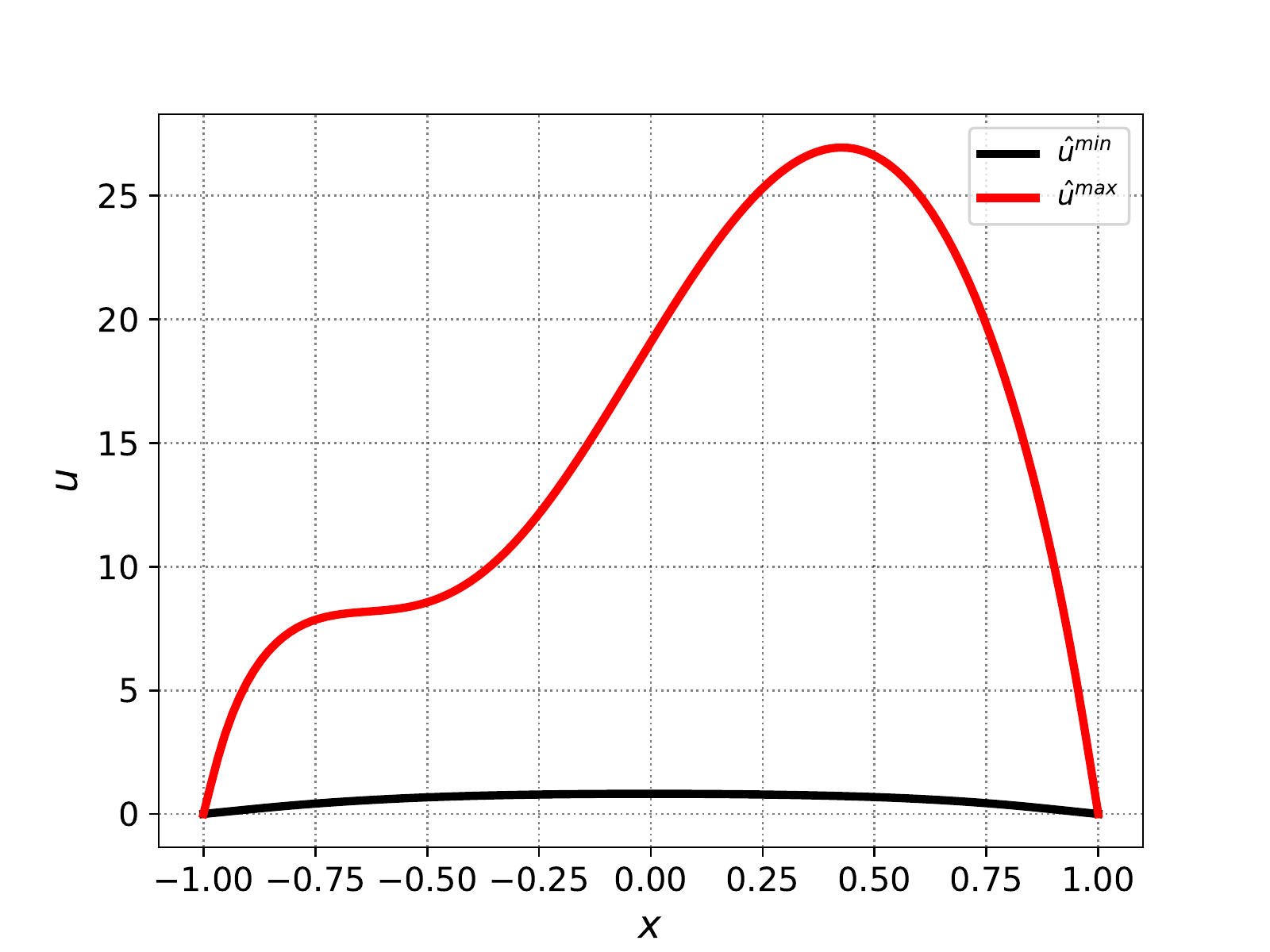}
\caption{$t=0.7s$}
\end{subfigure}

\begin{subfigure}{1.0\linewidth}
\centering
\includegraphics[scale=0.35]{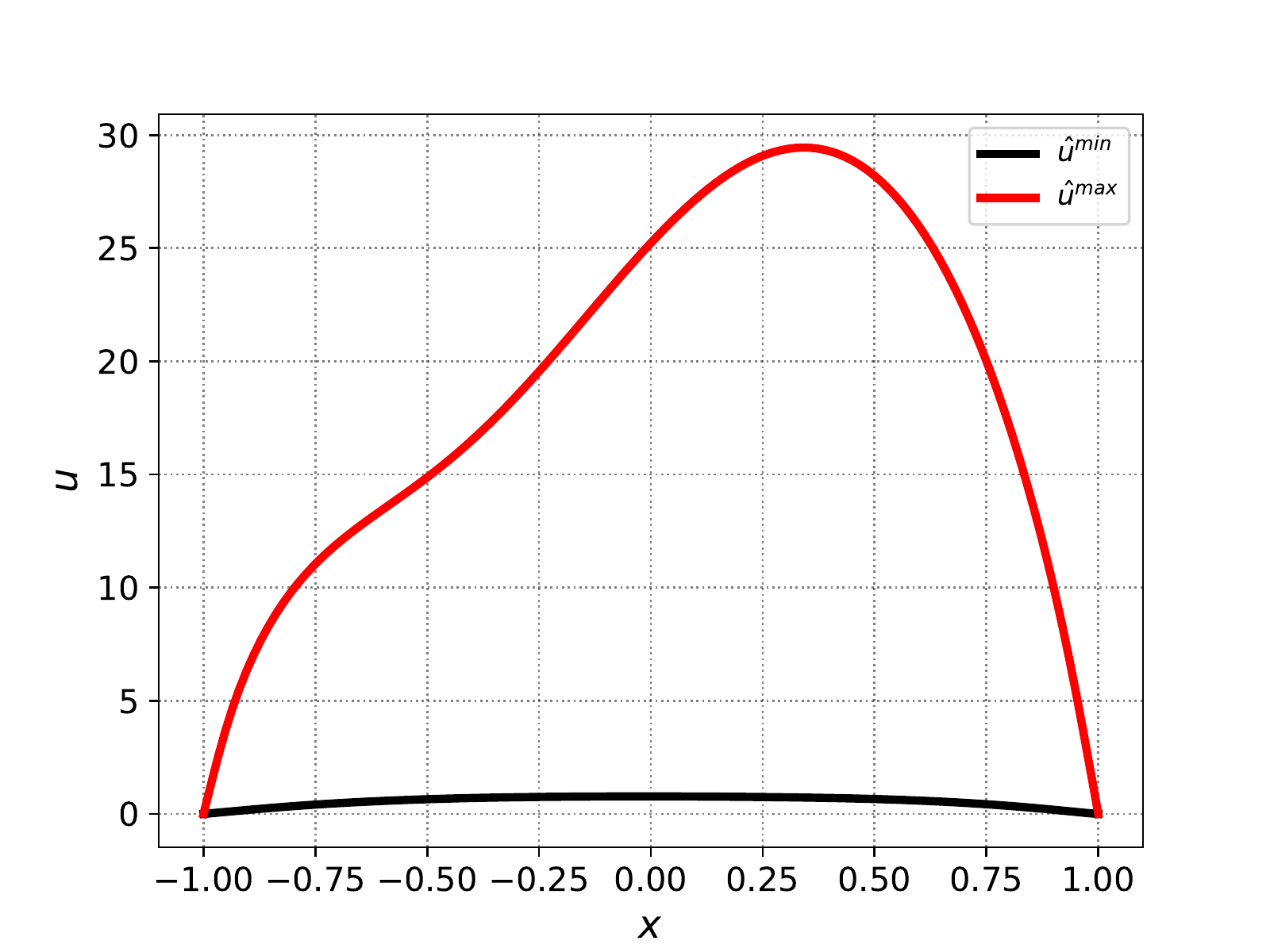}
\caption{$t=1.0s$}
\end{subfigure}
    \caption{Maximum and minimum predicted output fields for the nonlinear PDE problem for different points in time. Note the change in axis values.}
    \label{fig:nonPDEUfield}
\end{figure}

\begin{figure}
\begin{subfigure}{.5\linewidth}
\includegraphics[scale=0.35]{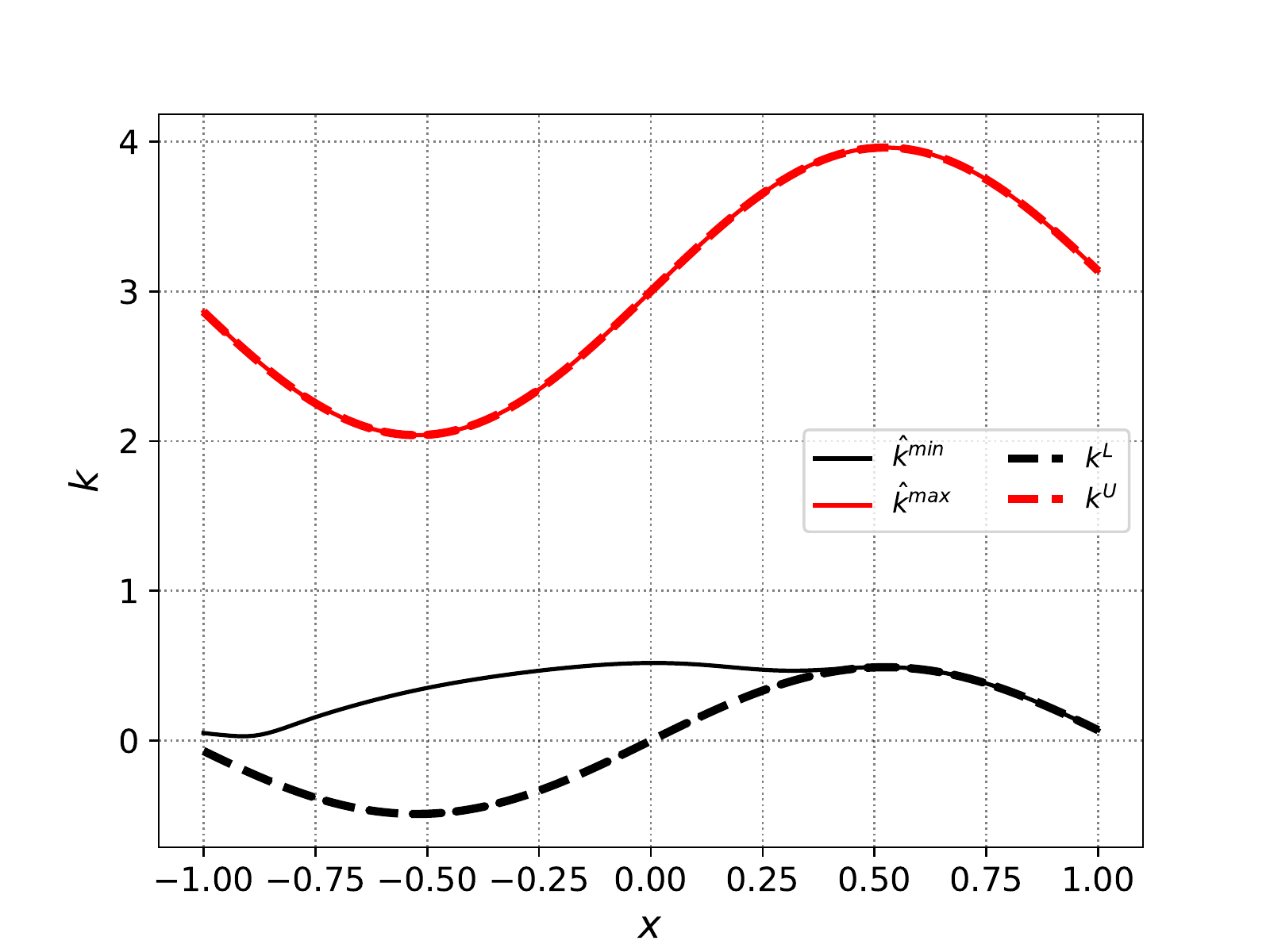}
\caption{$t=0.2s$}
\end{subfigure}
\begin{subfigure}{.5\linewidth}
\includegraphics[scale=0.35]{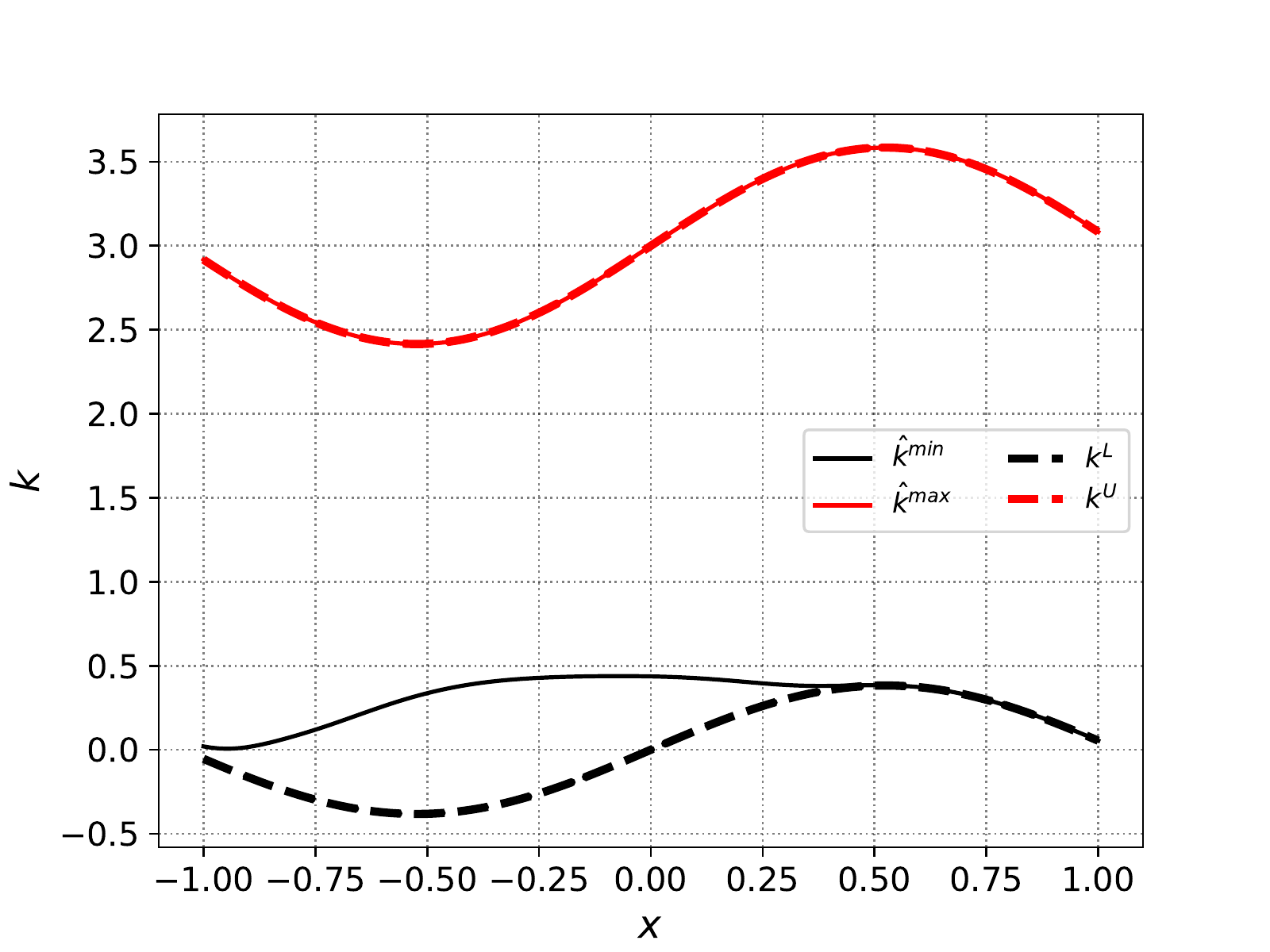}
\caption{$t=0.7s$}
\end{subfigure}
\begin{subfigure}{1.0\linewidth}
\centering
\includegraphics[scale=0.35]{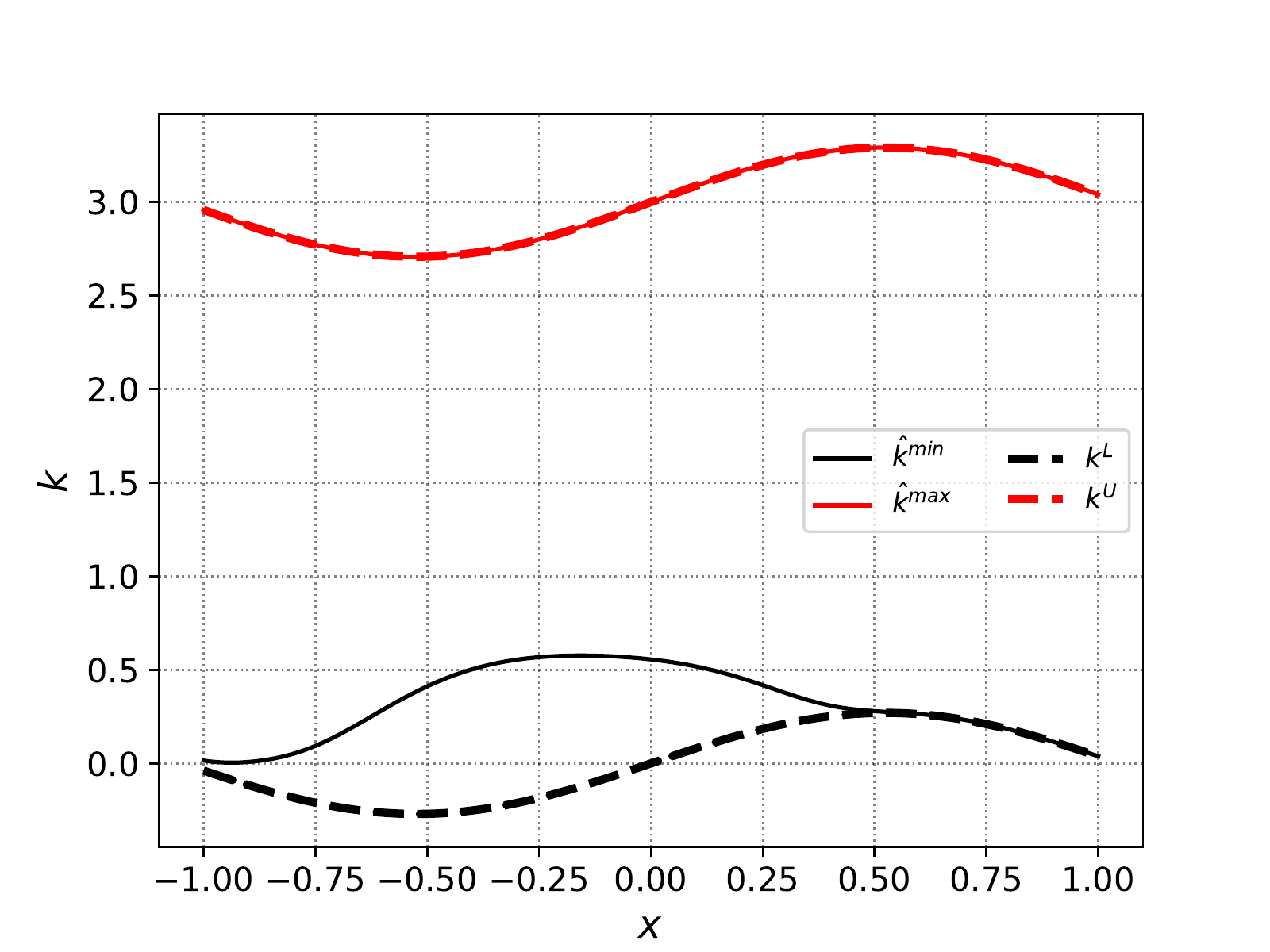}
\caption{$t=1.0s$}
\end{subfigure}
    \caption{Predicted and constraining input fields for the nonlinear PDE problem for specific time instances. Note the change in axis values.}
    \label{fig:nonPDEnoundinginputfield}
\end{figure}

\begin{figure}
\begin{subfigure}{0.5\linewidth}
\includegraphics[scale=0.35]{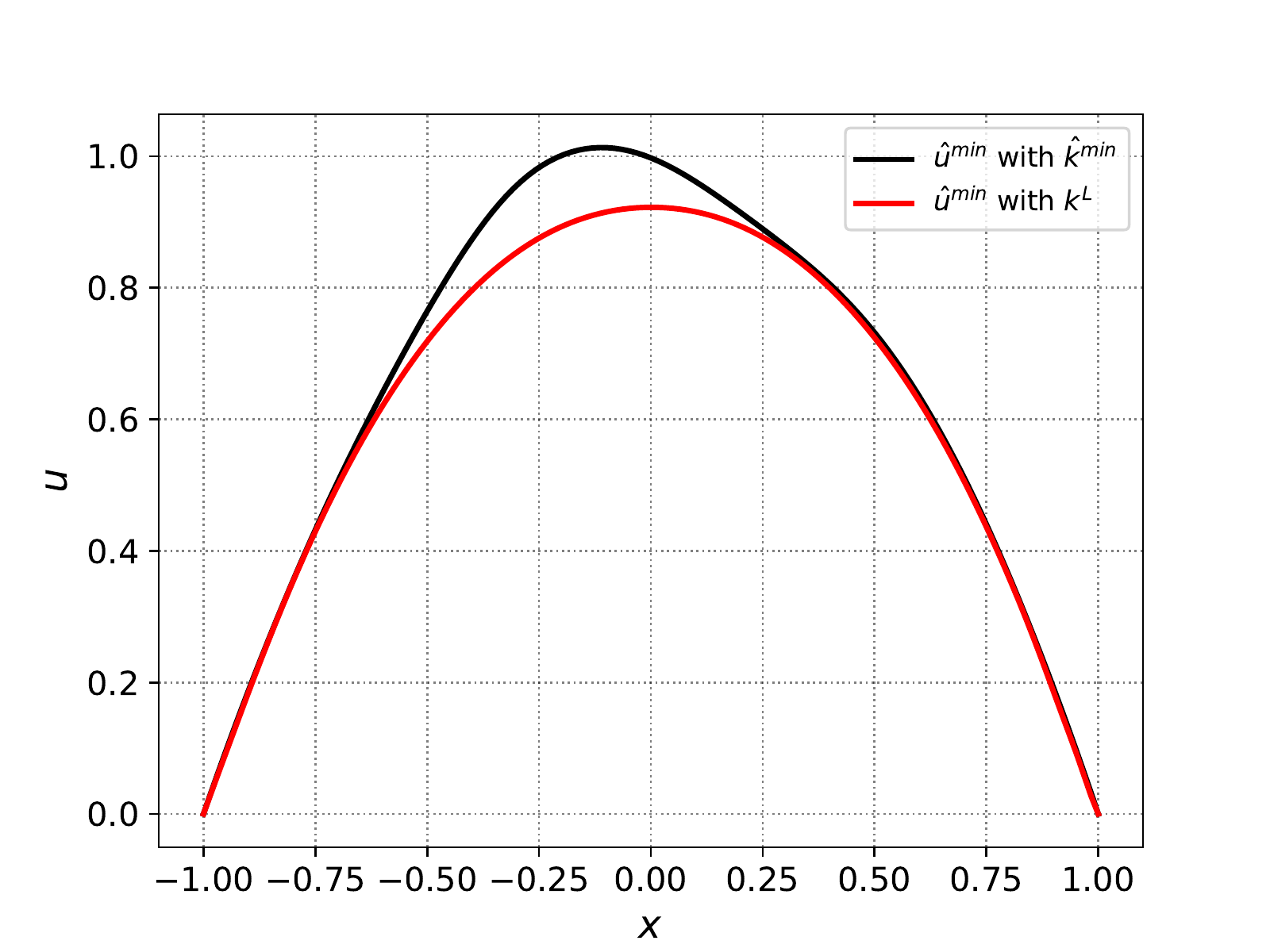}
\caption{$t=0.2s$}
\end{subfigure}
\begin{subfigure}{0.5\linewidth}
\includegraphics[scale=0.35]{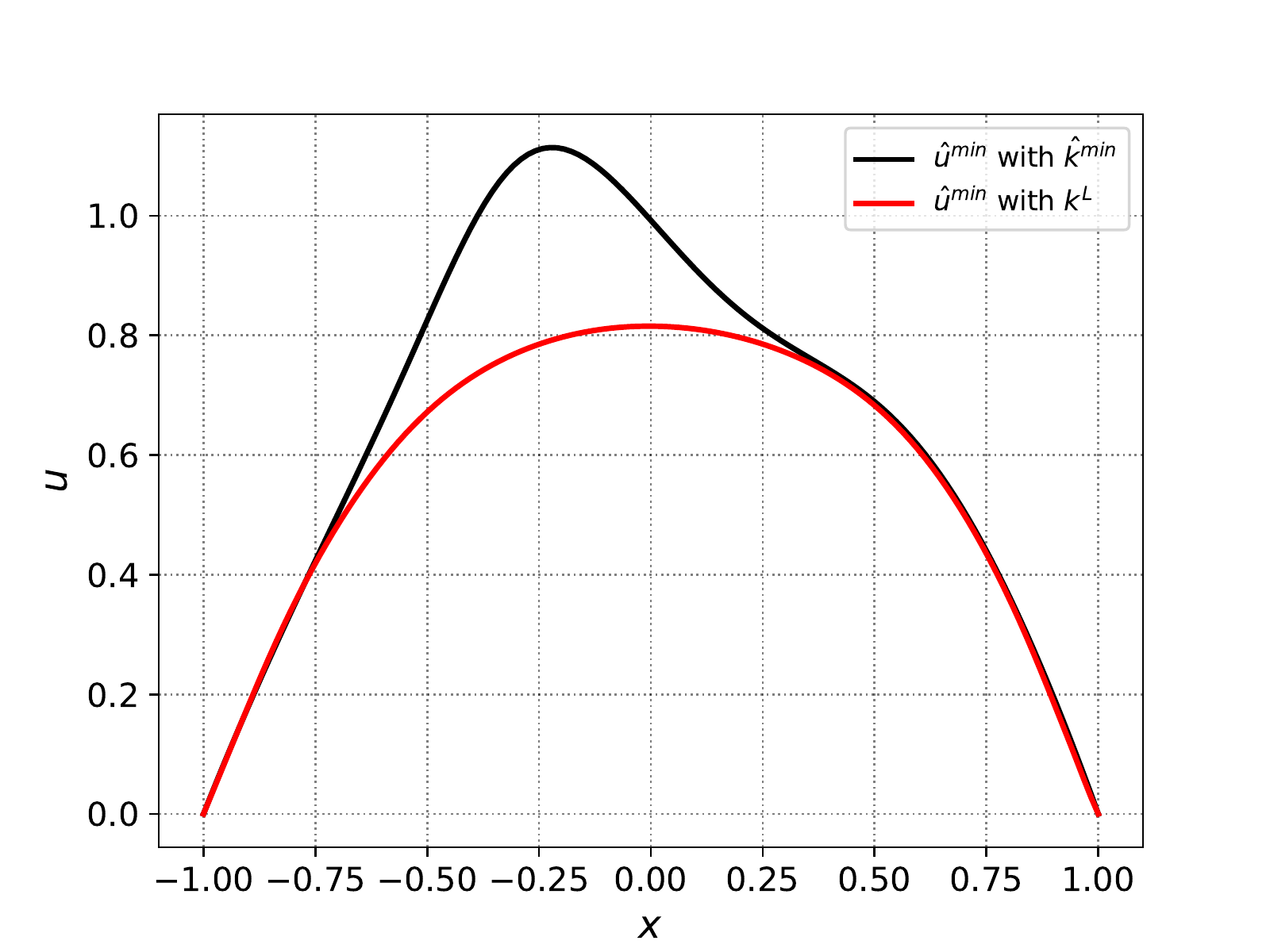}
\caption{$t=0.7s$}
\end{subfigure}
\begin{subfigure}{1.0\linewidth}
\centering
\includegraphics[scale=0.35]{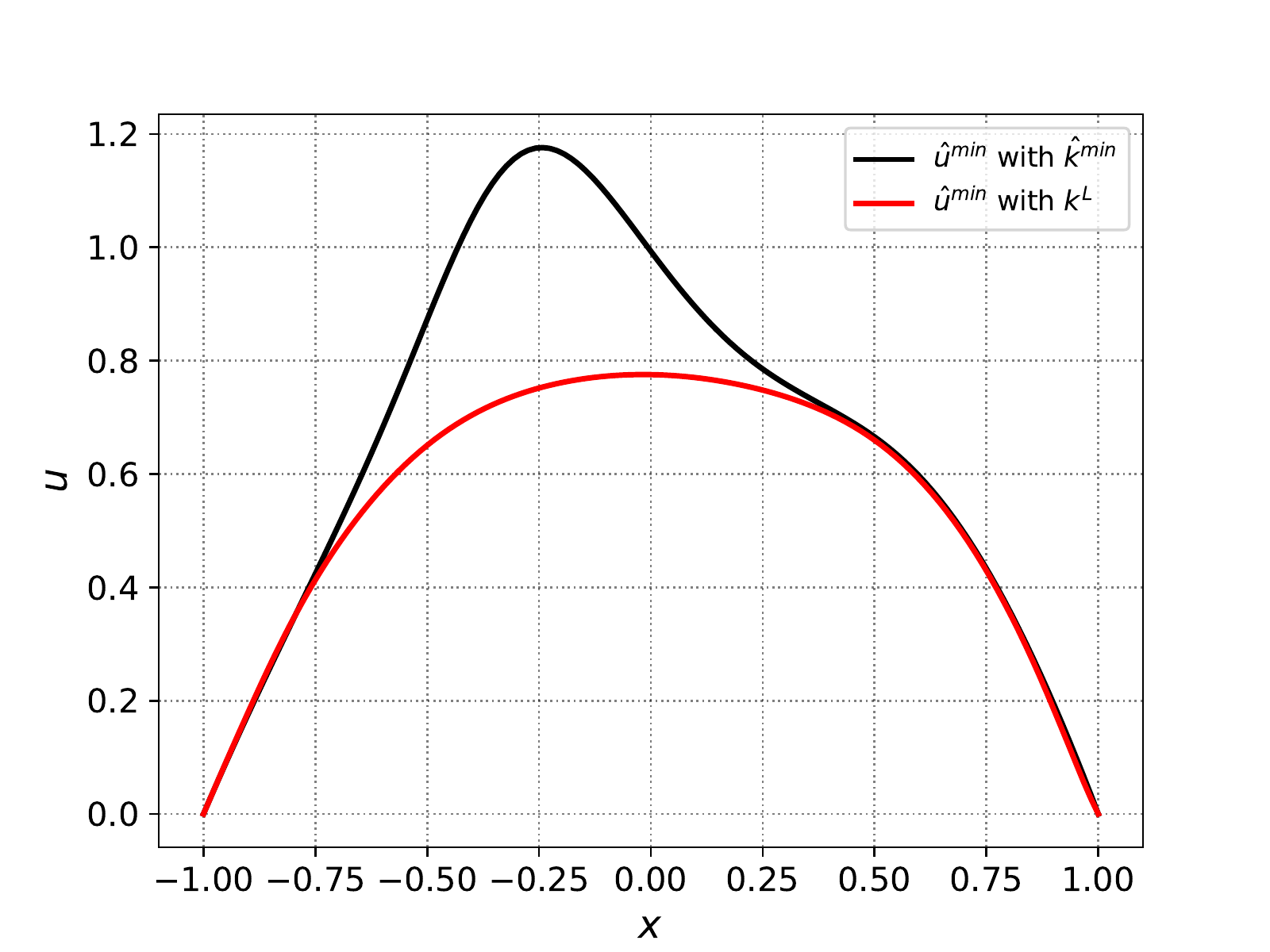}
\caption{$t=1.0s$}
\end{subfigure}
    \caption{Minimum predicted output field with corresponding minimum input field vs predicted output field using the minimum bound $k^{L}$ for nonlinear PDE problem.}
    \label{fig:nonPDECompuminubound}
\end{figure}

\begin{figure}
\begin{subfigure}{.5\linewidth}
\includegraphics[scale=0.27]{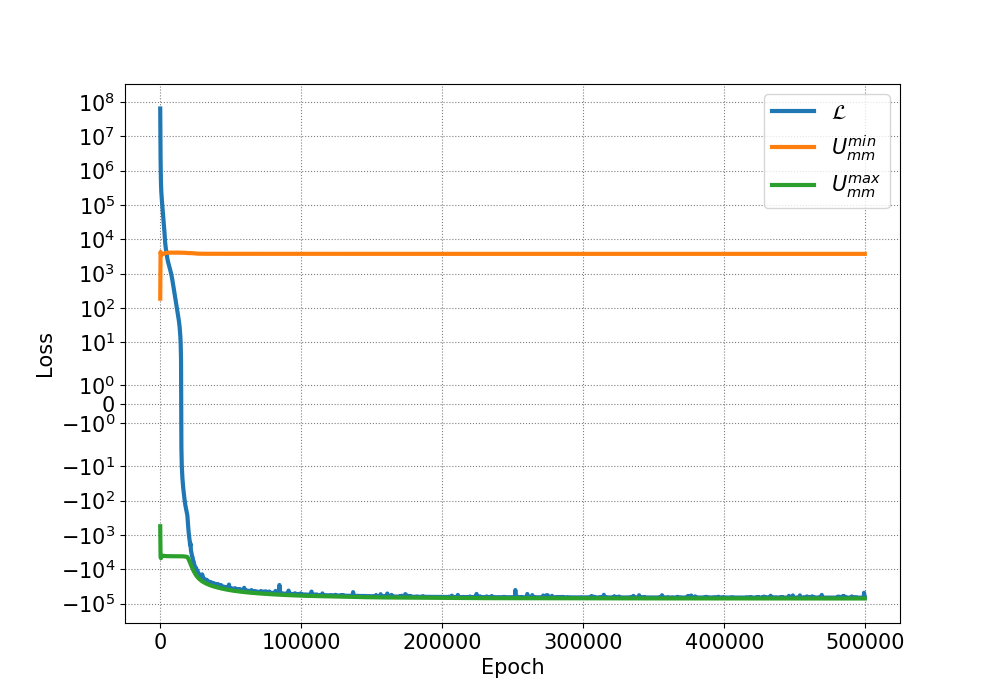}
\caption{}\label{fig:nonPDENonResErrors}
\end{subfigure}
\begin{subfigure}{.5\linewidth}
\includegraphics[scale=0.27]{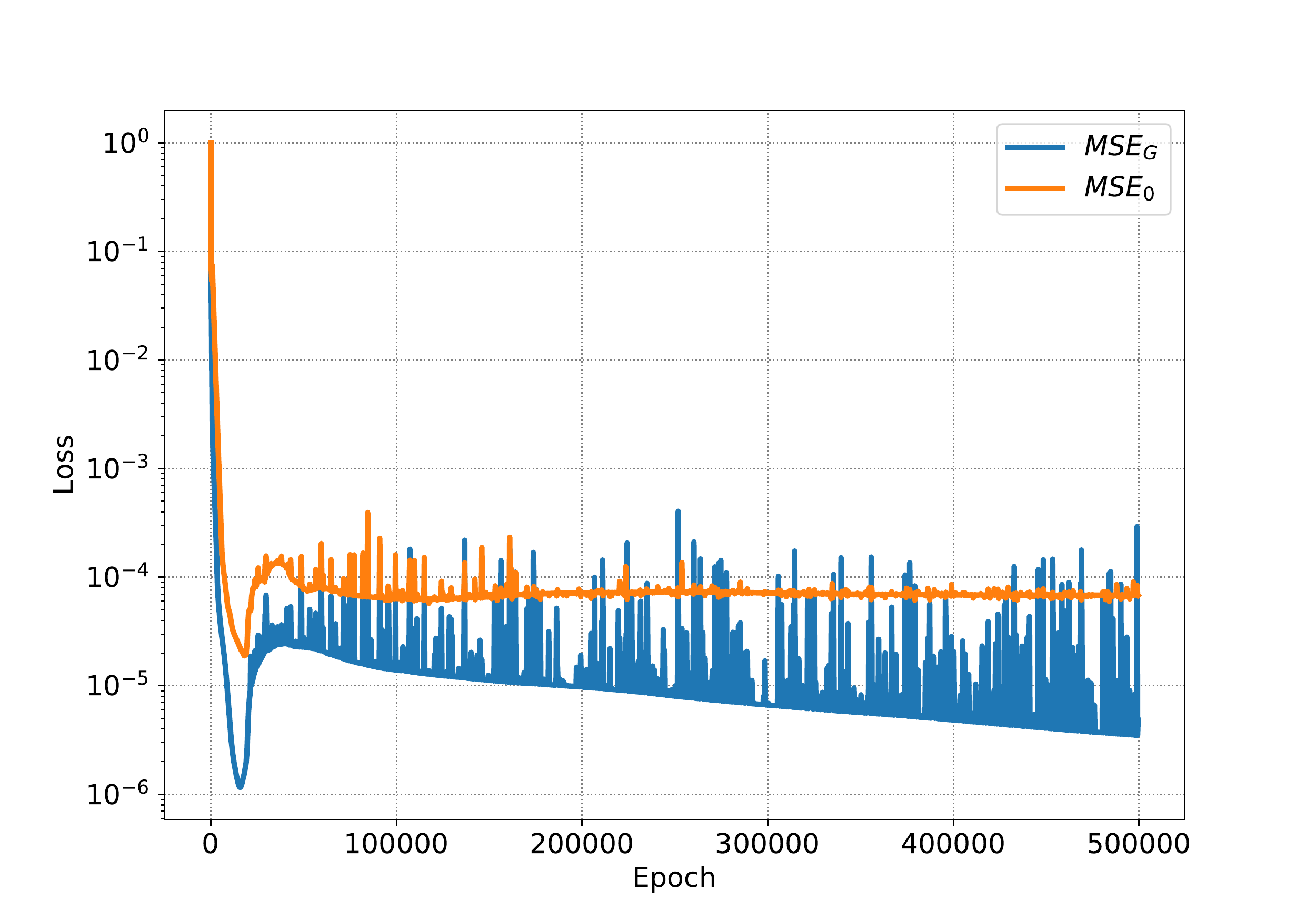}
\caption{}\label{fig:nonPDERESErrors}
\end{subfigure}
    \caption{(a) Non-residual losses over training epochs and (b) residual losses over training epochs for nonlinear interval PDE.}
    \label{fig:nonPDEErrors}
\end{figure}

The maximum and minimum bounding primary output fields are plotted for three different points in time in Figure \ref{fig:nonPDEUfield}. It can be seen that due to the interval parameter $k^{I}$ in equation (\ref{eq::nonPDE}), the primary output can take drastically different profiles. 
Upon solution, the corresponding bounding fields of the inputs for these three time points accompanied by the allowed upper and lower limits $k^{U}$ and $k^{L}$ are shown in Figure \ref{fig:nonPDEnoundinginputfield}. It can be noticed that $\hat{k}^{max}$, the approximated maximum field value, is approaching the upper limit $k^{U}$. On the other hand the approximated minimum field value $\hat{k}^{min}$ does not exactly correspond to the lower bound. In order to find this approximated field an appropriate sampling technique would have been required if one was to solve the problem using traditional approaches such as Monte-Carlo. With iPINN this field is directly available to us after solving the optimization problem of equation \ref{eq::WeightsWithLoss}.
Figure \ref{fig:nonPDECompuminubound} compares the primary output field solutions of eq. (\ref{eq::nonPDE}) when $k^{I}= \hat{k}^{min}$ and $k^{I}= \hat{k}^{L}$ to highlight that the lower bounding limit does indeed not yield the minimum primary output field.
Lastly, the non-residual and the residual errors of the problem over the training process are shown in Figure \ref{fig:nonPDEErrors}. It can be seen that after around $100,000$ epochs the non-residual values have converged to a stable value (Figure \ref{fig:nonPDENonResErrors}) where $U_{mm}^{min}$ and $U_{mm}^{max}$ correspond to the minimum and maximum field components of $U_{mm}$.
Furthermore, after $100,000$ epochs the 
initial $\text{MSE}_{0}$ and PDE residuals $\text{MSE}_{G}$ have both reached very accurate values lower than $1e-4$, see Figure \ref{fig:nonPDERESErrors}.

\section{Discussion and outlook}\label{sec::5}
The presented paper proposes a physics-informed neural network formulation for solving interval PDEs which can be used to obtain approximate solutions to fuzzy PDEs with convex membership functions as well. The proposed network infrastructure, termed interval physics-informed neural networks (iPINNs), consists of two separate feedforward neural networks with one aiming to approximate minimum and maximum possible primary outputs, and the other one aiming to obtain the accompanying input fields that lead to the  specific solutions of the PDE.
The presented approach is studied in an introductory example as well as in a one-dimensional structural problem and in a nonlinear time-dependent PDE, showing promising results. In particular the input fields that correspond to the bounds of the interval displacement field are obtained as a byproduct of the iPINN formulation without any necessary prior knowledge of the spatial correlation of the input field. We believe this framework has a lot of potential by being able to avoid the major problems of the finite element method when dealing with interval and fuzzy fields.

In the future we aim to extend this framework to higher dimensions and try to tackle inverse problems related to finding possible input intervals from observed output ranges. Furthermore we could potentially study more theoretical problems that are hard to study with interval finite element methods, i.e.  physical law discovery and higher dimensional nonlinear PDEs.

\bibliography{bib.bib}
\end{document}